\documentclass[11pt]{amsart}
\usepackage{amsbsy,amsfonts}
\usepackage{latexsym,eucal,amsmath,amsthm}
\usepackage{amssymb,psfig}

\newcommand{\db}{\mathbb }

\newcommand{\F}{\mbox{{\em F}}}

\newcommand{\Schwartz}{\mathcal S}

\def\R{{\db R}}

\def\F{{\mathcal C}}

\def\complex{\mathbb C}
\def\reals{\mathbb R}

\def\rp{{^{-1}}}
\def\p{\partial}

\def\eps{\varepsilon}
\def\be#1{\begin{equation} \label{#1}}
\def\bas{\begin{align*}}

\newtheorem{theorem}{Theorem}

\newtheorem{lemma}{Lemma}
\newtheorem{corollary}[theorem]{Corollary}

\theoremstyle{definition}

\theoremstyle{remark}
\newtheorem{remark}{Remark}

\numberwithin{equation}{section}
\numberwithin{lemma}{section}
\numberwithin{remark}{section}

\begin{document}

\title[Ill-posedness for NLS and NLW]
{Ill-posedness for nonlinear\\ Schr\"odinger and wave equations}

\author{Michael Christ}
\thanks{M.C.\ is supported in part by N.S.F. grant DMS 9970660.}
\address{University of California, Berkeley}
\author{James Colliander}
\thanks{J.C.\ is supported in part by N.S.F. grant DMS 0100595,
  N.S.E.R.C. grant RGPIN 250233-03 and a Sloan Fellowship.}
\address{\small University of Toronto}
\author{Terence Tao}
\thanks{T.T.\ is a Clay Prize Fellow and is supported in part by a grant
from the Packard Foundations.}
\address{University of California, Los Angeles}


\begin{abstract}
The nonlinear wave and Schr\"odinger equations on $\R^d$, 
with general power nonlinearity and with both the focusing and defocusing signs, 
are proved to be ill-posed in the Sobolev space $H^s$ whenever 
the exponent $s$ is lower than that predicted by scaling or Galilean invariances, or when the regularity is too low to support distributional solutions.  
This extends previous work \cite{cct1} of the authors, 
which treated the one-dimensional cubic nonlinear Schr\"odinger equation.  
In the defocusing case soliton or blowup examples are unavailable, and 
a proof of ill-posedness requires the construction of other solutions.
In \cite{cct1} this was achieved using certain long-time asymptotic behavior 
which occurs only for low power nonlinearities.
Here we analyze instead a class of solutions for which the zero-dispersion limit 
provides a good approximation.
The method is rather general and should 
be applicable to wider classes of nonlinear equations.
\end{abstract}

\date{August 13, 2003}

\maketitle

\section{Introduction}  \label{section:intro}

This paper is concerned with the low regularity behavior (and in particular ill-posedness) of the Cauchy problem for the generalized nonlinear Schr\"odinger equation 
\begin{equation} \label{gNLSp} \tag{gNLS}
\left\{
\begin{aligned}
-iu_t(t,x) + \Delta_x u(t,x) &= \omega |u|^{p-1}u(t,x)
\\
u(0,x)&=u_0(x)\in H^s(\reals^d)
\end{aligned}
\right.
\end{equation}
and the (complex) nonlinear wave equation
\begin{equation} \label{NLW} \tag{gNLW}
\left\{
\begin{aligned}
\square u(t,x) &= \omega |u|^{p-1} u(t,x)
\\ 
u(0,x) &= u_0(x)\in H^s(\reals^d)
\\
\p_t u(0,x) &= u_1(x)\in H^{s-1}(\reals^d)
\end{aligned}
\right.
\end{equation}
in $\reals^d$, where $\omega = \pm 1$ and $p > 1$, and $u: \reals \times \reals^d \to \complex$ is a complex-valued field.  Here $\Delta_x$ denotes the Laplacian $\Delta_x := \sum_j \frac{\p^2}{\p x_j^2}$, while $\Box := -\partial_t^2 + \Delta_x$ is the d'Alembertian.  The sign $\omega = +1$ is referred to as the \emph{defocusing} case, 
while the sign $\omega = -1$ is \emph{focusing}.

We say that the NLS equation \eqref{gNLSp} is \emph{locally well-posed in $H^s$} if for every $u_0 \in H^s$ there exist a time $T = T(\| u_0 \|_{H^s}) > 0$ and a (distributional) solution 
$u: [-T,T] \times \reals^d \to \complex$ to \eqref{gNLSp} which is in
the space $C^0 ([-T,T]; H^s_x)$, and such that the solution map $u_0 \mapsto
u$ is uniformly continuous\footnote{A reasonable alternative is to require
the solution mapping to be continuous but not necessarily uniformly continuous,
as is the case for certain equations, such as Burgers' equation, and the Korteweg-de Vries
Benjamin-Ono equations in the periodic case.}
from $H^s$ to $C^0 ([-T,T]; H^s_x)$.
Furthermore, there is an additional space $X$ in which $u$ lies, such
that $u$ is the unique solution to the Cauchy problem in $C^0 ([-T,T];
H^s_x)\cap X$; and 
$X\subset L^p_{t,x,\text{loc}}$,
so $|u|^{p-1} u$ is a well-defined spacetime distribution.   
A similar notion of local well-posedness can be formulated
for \eqref{NLW}, with the initial datum 
$(u_0,u_1)$ in $H^s_x \times H^{s-1}_x$, 
and $u$ itself in $C^0([-T,T]; H^s_x) \cap C^1([-T,T];H^{s-1}_x)$.
Here $H^s$ denotes the usual inhomogeneous Sobolev space.

The sharp range of exponents $p$, $s$, $d$ for which one has local well-posedness for \eqref{gNLSp} and \eqref{NLW} has been almost completely worked out.  For \eqref{gNLSp}, 
the scaling symmetry 
\be{nls-scaling}
u(t,x) \mapsto \lambda^{-2/(p-1)} u(\frac{t}{\lambda^2}, \frac{x}{\lambda})
\end{equation}
for $\lambda > 0$ leads one to the heuristic constraint
$$ s \geq s_c := \frac{d}{2} - \frac{2}{p-1}$$
(since the scaling leaves $\dot H^{s_c}$ invariant), while the Galilean invariance
\be{gal}
u(t,x) \mapsto e^{-iv \cdot x/2} e^{i |v|^2 t/4} u(t, x-vt)
\end{equation}
for arbitrary velocities $v \in \R^d$ leaves the $L^2$ norm invariant, and similarly leads to the heuristic
$$ s \geq 0.$$

Cazenave and Weissler \cite{cwI} showed that one indeed has local well-posedness when $s \geq 0$ and $s > s_c$, although in the case when $p$ is not an odd integer, we impose a natural compatibility condition $p > \lfloor s \rfloor + 1$ in order to ensure that the nonlinear term has sufficient smoothness.\footnote{
The assumption $p > \lfloor s \rfloor + 1$  is imposed for technical reasons
only; it might conceivably be relaxed or removed. See for instance
\cite{kato} for some work in this direction.} (In the case 
$s=s_c \geq 0$ one also has local well-posedness, but now the time of existence $T$ depends on the datum itself rather than merely on its $H^s$ norm.  See \cite{cwI}.)

Meanwhile, for \eqref{NLW}, the scaling symmetry
\be{wave-scale}
u(t,x) \mapsto \lambda^{-2/(p-1)} u(\frac{t}{\lambda}, \frac{x}{\lambda})
\end{equation}
once again gives the constraint
$$ s \geq s_c := \frac{d}{2} - \frac{2}{p-1}.$$
The equation $\square u = \omega |u|^{p-1}u$ also has the
Lorentz symmetries
\be{lorentz-scale}
u(t,x) \mapsto 
u\big(
\frac{t-vx_1}{(1-|v|^2)^{1/2}},
\frac{x_1-vt}{(1-|v|^2)^{1/2}},
x_2,\dots,x_d\big),
\end{equation}
for all subluminal speeds $-1< v < 1$. 
Certain combinations of these with scaling symmetries
give rise to the heuristic
$$ s \geq s_{{\rm conf}} := \frac{d+1}{4} - \frac{1}{p-1},$$
see e.g.\ \cite{sogge:wave}. 

Lindblad and Sogge \cite{lindbladsogge:semilinear} (see also
\cite{kapitanski:wayw},  \cite{tao:keel}) obtained local
well-posedness results for \eqref{NLW} 
analogous to those of \cite{cwI} for \eqref{gNLSp}:
\eqref{NLW} is well-posed in $H^s (\R^d), d \geq 2,$ whenever $s \geq
\max(s_c, s_{{\rm conf}})$ 
under three additional technical assumptions.
First, for smoothness reasons, one imposes the compatibility condition $p > \lfloor s \rfloor
+ 1$  when $p$ is not an odd integer.  Second, one needs to impose a
condition  on $s$ in order for the nonlinearity to make 
sense as a distribution; for $d \geq 2$, this condition is that $s
\geq 0$ and in dimension $d=1$, $s \geq \max(0,\frac{1}{2} -
\frac{1}{p})$.
Third, one also assumes an (apparently artificial) additional condition
$$p(\frac{d+1}{4} - s) \leq \frac{d+1}{2d} (\frac{d+3}{2} - s),$$
which only becomes relevant in dimensions $d \geq 4$ and for very low
values of $s$ and $p$.  
This artificial condition has been improved slightly \cite{tao:lowreg} to
$$p(\frac{d}{4} - s) \leq \frac{1}{2} (\frac{d+3}{2} - s).$$ 
Once again, at the endpoint $s=s_c \geq s_{{\rm conf}}$, 
local well-posedness holds, in the weaker sense that the time of existence is allowed to 
depend on the datum itself rather than merely on its norm.

The purpose of this paper is to complement these positive results
with ill-posedness results for most other Sobolev exponents $s$; we have satisfactory results in the case of NLS and partial progress in the case of NLW.
Such ill-posedness results are already known for the focusing
versions $\omega = -1$ of \eqref{gNLSp} and \eqref{NLW};  our arguments
apply also to defocusing equations and indeed are insensitive to the
distinction between the two.\footnote{We note however 
that Lebeau \cite{lebeau} has proved an interesting strong instability result, 
which implies ill-posedness for the defocusing wave equation, 
for energy space-supercritical cases in $\reals^3$.}
This thus extends our previous paper \cite{cct1}, which treated
the 1D cubic NLS equation $(p,d) = (3,1)$, as well as the related KdV and mKdV equations.

\subsection{Schr\"odinger ill-posedness results.}

We now discuss the known ill-posedness results for the Schr\"odinger equation in the focusing case $\omega = -1$.  There are three cases: the $L^2$-subcritical case 
$s_c < 0$ ($p < 1 + \frac{4}{d}$), 
the $L^2$-critical case $s_c = 0$ ($p = 1 + \frac{4}{d}$), 
and the $L^2$-supercritical case $s_c > 0$ ($p > 1 + \frac{4}{d}$).

In the focusing $L^2$-supercritical case $s_c > 0$,
blowup in finite time from smooth data can be established via the virial identity (see
\cite{VPT}, \cite{Z72}, \cite{Glassey} or the textbook
\cite{SulemSulem}). 
By combining this blowup example with the scaling \eqref{nls-scaling} one can show 
blowup in arbitrarily short
time for data in $H^s$ when $s < s_c$, thus complementing the results
in \cite{cwI}. 

In the focusing $L^2$-critical case $s_c = 0$, a similar blowup example can be produced
by taking a ground state solution and applying the pseudo-conformal transformation, and by rescaling this example one obtains blowup in arbitrarily short time when $s < s_c = 0$.

In the focusing $L^2$-subcritical case $s_c < 0$, one has global 
well-posedness in $L^2$, and so there are no examples of blowup solutions 
with which to repeat the previous scaling argument.  Nevertheless, manipulation of
soliton solutions via scaling and Galilean transformations demonstrates
that the solution map is not uniformly continuous in $H^s$ when $s < 0$, thus 
establishing a weaker form of ill-posedness in these spaces.
See the papers of Birnir et.\ al.\ \cite{birnir}, Kenig-Ponce-Vega \cite{kpvill},
and Biagioni-Linares \cite{biagionilinares} for this type of argument.
In the defocusing case, one lacks any useful explicit solutions on which to base
this type of reasoning.

Thus in the focusing case $\omega = -1$, there is
a satisfactory collection of ill-posedness examples which show that
the well-posedness results in \cite{cwI} are sharp.  Our first main
result extends these ill-posedness results to the defocusing case
$\omega = +1$, for which blowup solutions\footnote{It is known that
  there are no blowup solutions in the defocusing case, in the
  $H^1$-subcritical  case $s_c < 1$ (see e.g.\ \cite{cwI}).  The
  situation  when $s_c \geq 1$ remains open at present; the question
  of blowup in the defocusing $H^1$-supercritical case $s_c > 1$ seems
  particularly intractable.} 
and soliton solutions are unavailable.
Throughout this paper
we use $C \gg 1$ and $0 < c \ll 1$ to denote various large 
and small constants respectively, which may vary from line to
line. Let $\Schwartz$ denote the Schwartz class.

\begin{theorem}\label{thm:NLS}
Let $p >1$ be an odd integer, let $d\ge 1$, and let $\omega = \pm 1$.
For any $s<\max(0,s_c)$, the Cauchy problem \eqref{gNLSp} fails to be well-posed in $H^s(\reals^d)$.
More precisely, for any $0< \delta, \epsilon <1$ and for any $t>0$
there exist solutions $u_1, u_2$ of \eqref{gNLSp} with initial data
$u_{1}(0), u_{2}(0) \in \Schwartz$ such that
\begin{equation}
  \label{datasize}
  {{\| u_{1}(0) \|}_{H^s}}, {{\| u_{2}(0) \|}_{H^s}} \leq C \epsilon,
\end{equation}
\begin{equation}
  \label{dataclose}
  {{\| u_{1}(0) - u_{2}(0) \|}_{H^s}} < C \delta,
\end{equation}
\begin{equation}
  \label{bang}
  {{\| u_1 (t) - u_2 (t) \|}_{H^s}} > c \epsilon.
\end{equation}
Thus the solution operator fails to be uniformly continuous on $H^s$.

If $p > 1$ is not an odd integer, then the same conclusion holds
provided that there exists an integer $k>d/2$ such that $p\ge k+1$ and
$s<\max(0,s_c)$.
\end{theorem}

\begin{figure}[htbp] \centering
 \ \psfig{figure=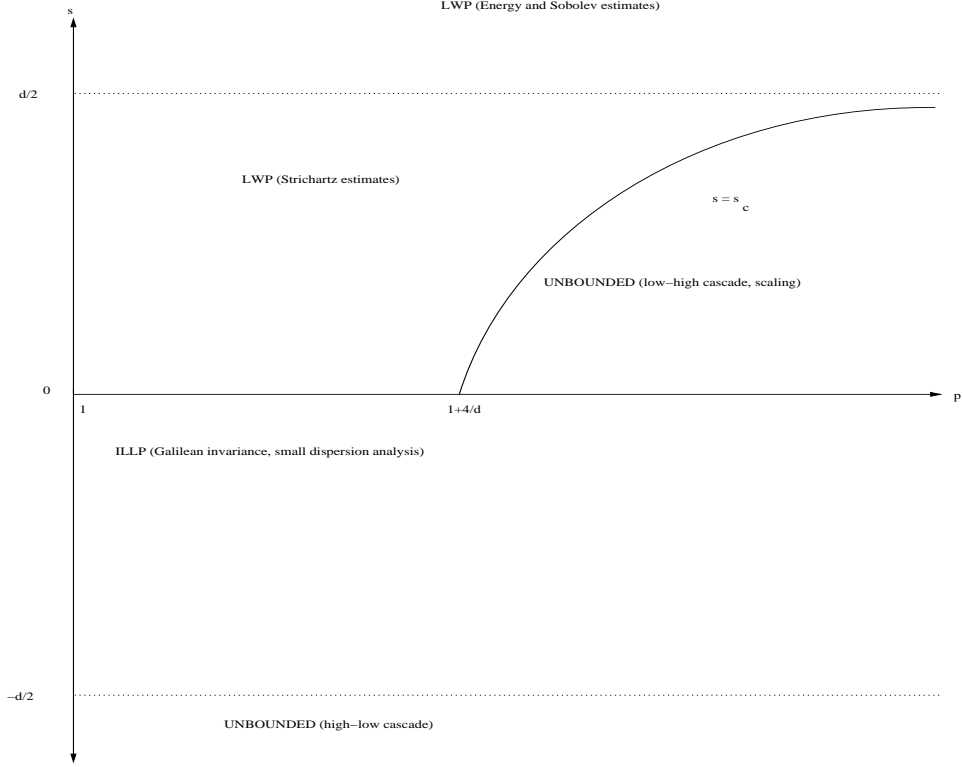,width=5in,height=4in}
 \caption{
A schematic depiction of the NLS results for a fixed $d$ and variable
$p$ and $s$; the diagram is not drawn to scale, and omits some
technicalities such as the presence of the integer $k$ when $p$ is not
an odd integer.  When $s \geq \max(s_c, 0)$ one has local
well-posedness by Strichartz estimates (or by energy methods when
$s>d/2$); see \cite{cwI}.  Theorem \ref{thm:NLS} uses primarily
Galilean invariance and small-dispersion analysis to yield
ill-posedness for $s < 0$.  When $s < -d/2$ or $0 < s < s_c$ we have
the stronger result (Theorem \ref{thm:bignorm}) that the $H^s$ norm
can become unbounded from arbitrarily small initial data in
arbitrarily small time; the two cases exploit a high-to-low and
low-to-high cascade of frequencies respectively in the small
dispersion analysis, as well as some use of the scale invariance. In
this figure and those that follow, we have suppressed the smoothness
restrictions like $p \geq \lfloor{s} \rfloor +1.$}
 \end{figure}

The presence of the integer $k$ is perhaps an artifact of our argument, 
which relies on the energy method and hence
requires the nonlinear term $|u|^{p-1} u$ to be sufficiently regular. 

In \cite{cct1} we obtained this result in the special case $(p,d) = (3,1)$ of the 1D cubic NLS.  Our arguments there relied mainly on the existence of modified scattering solutions 
for that equation (and on their failure to behave like unmodified scattering solutions),
and so in particular do not apply to the case $d=1$, $p>3$.  
Our methods here are different (and somewhat simpler), relying instead on some 
quantitative analysis of the NLS equation
\begin{equation} \label{smalldispersion}
\left\{
\begin{aligned}
-i\phi_s(s,y) + \nu^2 \Delta_y \phi(s,y) &= \omega |\phi|^{p-1} \phi(s,y)
\\
\phi(0,y)&=\phi_0(y)
\end{aligned}
\right.
\end{equation}
for fixed initial datum $\phi_0$, in the \emph{small dispersion} regime $\nu \to 0$; 
\eqref{smalldispersion} can be transformed back into \eqref{gNLSp} via
a suitable rescaling of space and time, namely 
\be{u-phi}
u(t,x) := \phi(t, \nu x).
\end{equation}
Combining this transformation with \eqref{nls-scaling} and \eqref{gal} 
yields a family of solutions with three independent parameters $\nu,\lambda,v$. 
The small dispersion initial datum $\phi_0$ in our analysis varies only over the 
one-parameter family of scalar multiples of a fixed function.
Various cases of our results will be proved by choosing different combinations
of parameter values.

In \cite{cct1}, small dispersion analysis
was used to handle the supercritical case $s < s_c$; we reproduce that argument here and 
show how it can also handle the super-Galilean case $s < 0$ via a
frequency modulation argument related to the Galilean transform \eqref{gal}.  
We prove Theorem~\ref{thm:NLS} in Sections~\ref{section:negatives} and
\ref{section:positives}, after some 
small dispersion analysis in Section~\ref{section-zero}.

Koch and Tzvetkov \cite{kochtzvetkovill}
have shown by an argument exploiting small dispersion considerations
that the solution map for the Benjamin-Ono
equation fails to be uniformly continuous in $H^s(\reals)$ for $s>0$, although it
is continuous for sufficiently large $s$.

The form of ill-posedness demonstrated in Theorem~\ref{thm:NLS} is rather mild; 
the solution map is merely proved not to be uniformly continuous. 
The next result establishes much worse behavior,
involving rapid growth of norms, in certain cases. 
We will call this type of behavior {\em norm inflation}.

\begin{theorem}  \label{thm:bignorm} Let $p > 1$ be an odd integer, let $d\ge 1$, and let $\omega = \pm 1$. 
Suppose that either $0<s< s_c = \tfrac{d}2 -\tfrac{2}{p-1}$ or $s \le -d/2$.
Then for any $\eps>0$ there exist a solution $u$ of \eqref{gNLSp}
and $t\in\reals^+$ such that
$u(0)\in\Schwartz$,
\begin{align}
&\|u(0)\|_{H^s}< \eps,
\\
&0<t<\eps,
\\
&\|u(t)\|_{H^s}> \eps\rp.
\end{align}
In particular, for any $t>0$
the solution map $\Schwartz\owns u(0)\mapsto u(t)$
for the Cauchy problem \eqref{gNLSp} 
fails to be continuous at $0$ in the $H^s$ topology. 

If $p\ge 1$ is not an odd integer, then the same conclusion holds
provided that there exists an integer $k>d/2$ such that $p\ge k+1$ and 
$-k<s<\max(0,s_c)$.
\end{theorem}

The mechanism underlying norm inflation is the transfer of energy
from low to high Fourier modes for $s>0$, and from high to low modes
for $s\le -\frac{d}{2}$.  We prove this theorem in Sections \ref{section:negatives}
and \ref{section:verylows}, using the small dispersion analysis from 
Section~\ref{section-zero}.

A related result, for nonlinear wave equations, was obtained a number of years
ago by Kuksin \cite{kuksin} via a closely related argument, involving
approximation of the PDE by an ODE in a suitable small dispersion/large time regime.

\subsection{Wave ill-posedness results.}

We now consider ill-posedness of the nonlinear wave equation \eqref{NLW}.
In the focusing case $\omega = -1$, it is well-known that blowup and 
ill-posedness can be obtained via the ODE method.  Indeed, the blowup of solutions of
the ODE $u_{tt} = +|u|^{p-1} u$ means that there are constant-in-space 
solutions to \eqref{NLW} which blow up in finite time.  
Truncating the initial data in space yields compactly supported solutions which blow up in 
finite time, by virtue of the finite speed of propagation. 
Transforming these blowup solutions using either the scaling symmetry
\eqref{wave-scale} or the Lorentz symmetry \eqref{lorentz-scale} 
establishes blowup in arbitrarily short time when $s < s_c$ or $s < s_{{\rm conf}}$ 
respectively. See \cite{sogge:wave}, \cite{lindbladsogge:semilinear}. 

In the defocusing case $\omega = +1$, the blowup solutions are unavailable 
(at least when $s_c \leq 1$; the question of blowup in the supercritical case $s_c > 1$ 
remains open).  Nevertheless, a small dispersion analysis 
can still yield ill-posedness as in Theorem \ref{thm:NLS}.  
For simplicity we treat only the defocusing case.

\begin{theorem}  \label{thm:NLW}
Let $d\ge 1$, $\omega = + 1$,  and $p > 1$.
If $p$ is not an odd integer, we assume in addition that $p \geq k+1$
for some integer $k>d/2$.
Then the Cauchy problem \eqref{NLW} is ill-posed in $H^s(\reals^d)$
for all $ s<s_c = \tfrac{d}2-\tfrac2{p-1}$.
More precisely, for any $0< \delta, \epsilon <1$ and for any $t>0$
there exist solutions $u_1, u_2$ of \eqref{NLW} with initial data
$(u_{1}(0), u_{1t} (0)), (u_{2}(0), u_{2t}(0)) \in \Schwartz \times \Schwartz$ such that
\begin{equation}
  \label{nlwdatasize}
  {{\| u_{1}(0) \|}_{H^s}}, {{\| u_{2}(0) \|}_{H^s}} \leq  \epsilon, ~u_{1t} (0) = u_{2t}(0) = 0,
\end{equation}
\begin{equation}
  \label{nlwdataclose}
  {{\| u_{1}(0) - u_{2}(0) \|}_{H^s}} < C \delta,
\end{equation}
\begin{equation}
  \label{nlwbang}
  {{\| u_1 (t) - u_2 (t) \|}_{H^s}} > c \epsilon.
\end{equation}
In particular, for any $T>0$ and any ball $B\subset H^s\times H^{s-1}$,
it is not the case that the map from initial data to solution is uniformly continuous
as a map from $B$ to $C^0([-T,T],H^s(\reals^d))$.
\end{theorem}
Theorem~\ref{thm:bignorm-wave} and Corollary~\ref{1dcorollary}
below together imply a stronger form of illposedness for all $s<s_c$, {\em except} $s=0$.

\begin{figure}[htbp] \centering
 \ \psfig{figure=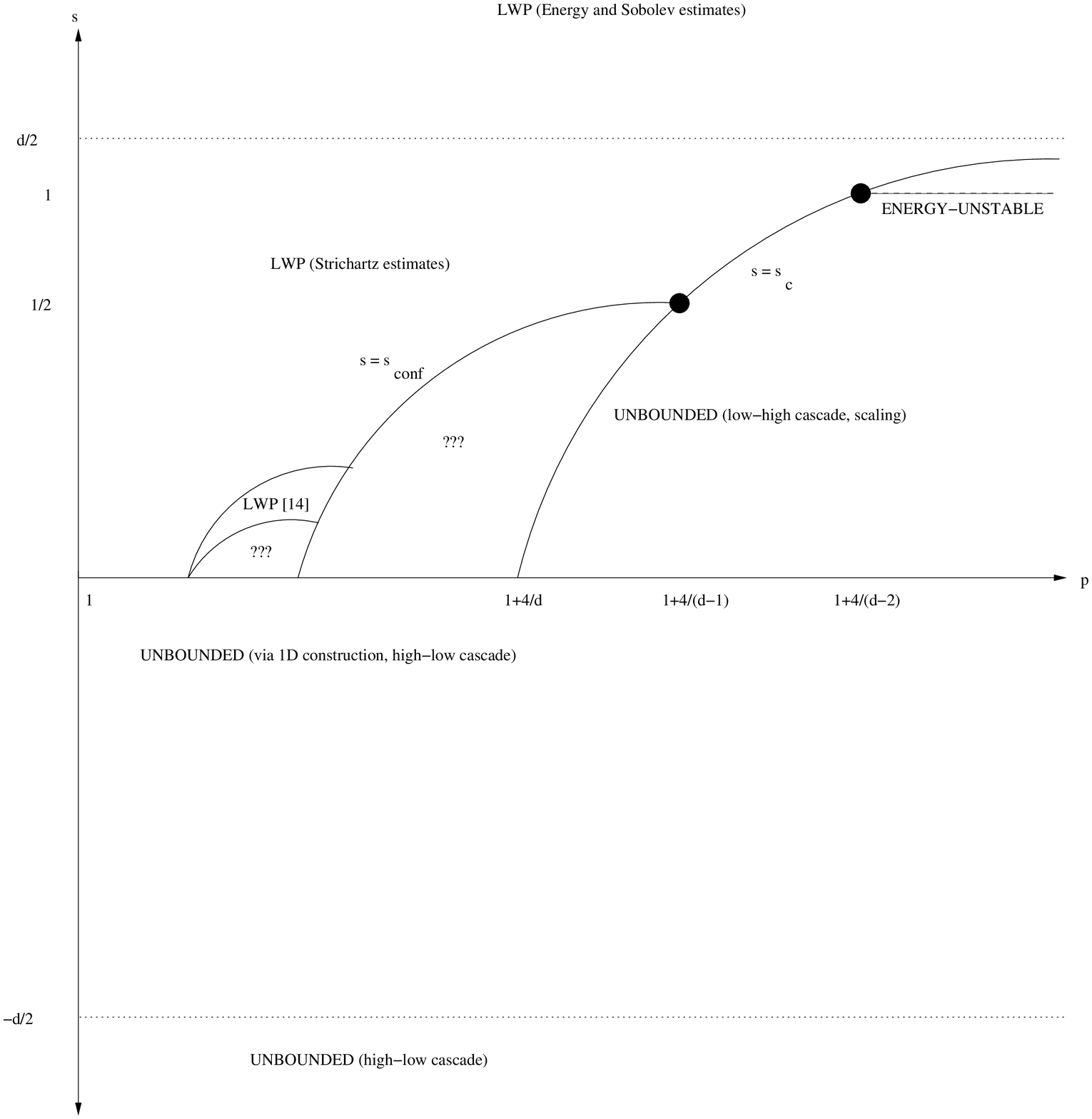,height=4in,width=5in}
 \caption{A schematic depiction of the defocusing
NLW results for $d > 1$.  
For $s < -d/2$ or $0 < s < s_c$ we have rapid norm inflation by Theorem \ref{thm:bignorm-wave},
and likewise for $s<0$ by Corollary~\ref{1dcorollary}.
For $\max(s_c,0) \leq s < s_{\rm conf}$ 
there is rapid norm inflation for $s<\tfrac12-\tfrac1p$ (not shown in the figure),
but the rest of this region remains open, though ill-posedness is known in the focusing case. 
The problem is illposed, with rapid decoherence, for $s=0$ by Theorem~\ref{thm:NLW}.
 In the energy supercritical cases $p > 1 + 4/(d-2)$ there is uniform boundedness, but
instability in the energy norm (\cite{brenner}, \cite{lebeau}, Theorem \ref{thm:energy-bad}).}
 \end{figure}

The analogue of Theorem \ref{thm:bignorm}, concerning norm inflation, also holds for \eqref{NLW}:

\begin{theorem}\label{thm:bignorm-wave} 
Let $d\ge 1$, $\omega = + 1$,  and $p >1$.
If $p$ is not an odd integer, we assume in addition that $p \geq k+1$
for some integer $k>d/2$.
Suppose that $0<s< s_c = \tfrac{d}2 -\tfrac{2}{p-1}$ or $s \le  -d/2$.
Then for any $\eps>0$ there exist a real-valued solution $u$ of \eqref{NLW}
and $t\in\reals^+$ such that $u(0)\in\Schwartz$,
\begin{align}
&\|u(0)\|_{H^s} < \eps,\\
& u_t(0) = 0
\\
&0<t<\eps,
\\
&\|u(t)\|_{H^s} > \eps\rp.
\end{align}
In particular, for any $t>0$
the solution map $\Schwartz\owns (u(0), u_t(0)) \mapsto (u(t), u_t(t))$
for the Cauchy problem \eqref{NLW} 
fails to be continuous at $0$ in the $H^s \times H^{s-1}$ topology.
\end{theorem}

In certain $H^1$-supercritical cases $0 \leq s \leq 1 < s_c$,
instability results were proven in \cite{brenner}, \cite{lebeau}.  
Our methods yield alternative proofs of a stronger result than that obtained in \cite{brenner},
and a weaker one than \cite{lebeau}.
For instance, suppose that $s_c > 1$, and consider the defocusing case $\omega = +1$.  
Because of the conserved energy
\be{energy-scale}
\int \tfrac{1}{2} |\nabla u|^2 + \tfrac{1}{2} |u_t|^2 + \tfrac{1}{p+1} |u|^{p+1}
\end{equation}
it is natural to consider well-posedness in the norm
$$\| (u, u_t) \|_X :=  \| \nabla u \|_2 + \| u_t \|_2 + \| u \|_{p+1}.$$
The energy conservation law thus means that the solution map is (formally) a bounded map 
from $X$ to itself.  Nevertheless, it cannot obey any sort of uniform continuity 
properties when $s_c > 1$:

\begin{theorem}  \label{thm:energy-bad}
Let $d\ge 3$, $\omega = + 1$, $p > 1$, and let $k > d/2$ be an integer.  If $p$ is not an odd integer, we impose the additional condition $p \geq k+1$.  Assume also that $s_c > 1$.  Then for any $0 < \eps, \delta < 1$ there exist $0<T<\eps$ and
real-valued solutions $u$, $u'$ defined on $[0,T]$ satisfying
\begin{align}
&\big\|(u(0), u_t(0))\big\|_X <C \eps,
\\
&\big\|(u(0), u_t(0)) - (u'(0), u'_t(0))\big\|_X <C \delta,
\\
&\big\|(u(T), u_t(T)) - (u'(T), u'_t(T))\big\|_X >c \eps.
\end{align}
\end{theorem}

We emphasize that our methods do not establish blowup for the $H^1$-supercritical defocusing wave equation, and the question of whether, say, smooth solutions stay smooth globally in time remains an interesting open problem.  However, the lack of continuity of the solution map in the energy class does suggest that the standard techniques used to obtain regularity will not be effective for this problem.

In the supercritical cases $s < s_c$ the proofs of Theorems 
\ref{thm:NLW}, \ref{thm:bignorm-wave}, and \ref{thm:energy-bad} proceed in analogy with their NLS counterparts, using 
small dispersion analysis and scaling arguments.  
However, we have not succeeded in adapting the decoherence argument in
\S\ref{section:negatives} to the regime $s_c<s<s_{\rm conf}$
for \eqref{NLW}. 

We now turn to NLW in the one-dimensional case $d=1$, which has some special features.  
In one dimension there are no local smoothing properties (because there is no decay in the fundamental solution), and so in particular there are no Strichartz estimates.  Thus, the only estimates available are those arising from energy methods and Sobolev embedding.  
This suggests that for $d=1$ there should be an additional restriction
$s\ge \frac{1}{2} - \frac{1}{p}$ 
for local well-posedness; this supersedes the
requirement $s\geq 0$ when $p > 2$.  
We thus introduce the threshold exponent
\begin{equation}
s_{\rm sob} = \max(0,\tfrac{1}{2} - \tfrac{1}{p}).  
\end{equation}
Its significance is that 
$H^{s_{\rm sob}}$ is the minimal 
Sobolev space for which Sobolev embedding places the solution of the 
linear wave equation in 
$C^0_TL^p_x$, so that the nonlinearity becomes locally integrable in the spatial variable.  
It is easy to check using energy methods (see also Section \ref{1d:sec}) that
NLW is locally well-posed when $s \ge s_{\rm sob}$; 
of course the usual assumption $p > \lfloor s_{\rm sob} \rfloor + 1$ 
must also be imposed
when $p$ is not an odd integer.    
Our final result confirms ill-posedness
for $s < s_{\rm sob}$, and establishes norm inflation
in the spirit of Theorem~\ref{thm:bignorm} for $s_c<s<s_{\rm sob}$.
Note that 
$\max (s_c, s_{\rm{conf}} )$ is always strictly less than $s_{\text{sob}}$.

\begin{theorem}\label{thm:bignorm-1d} 
Let $d = 1$, $\omega = \pm 1$.
Suppose either that $p$ is an odd positive integer, or that
$p > \lfloor s \rfloor + 1$.
Then \eqref{NLW} is ill-posed in $H^s$ for all 
$s< s_{\rm sob} = \max(0,\tfrac{1}2 -\tfrac{1}{p})$. 

More specifically,
if $s_c<s<s_{\rm sob}$ where $s_c = \tfrac12-\tfrac2{p-1}$,
then for any $\eps>0$ there exist a real-valued solution $u$ of \eqref{NLW}
and $t\in\reals^+$ such that $u(0)\in\Schwartz$,
\begin{align}
&\|u(0)\|_{H^s} < \eps,\\
& u_t(0) = 0
\\
&0<t<\eps,
\\
&\|u(t)\|_{H^s} > \eps\rp.
\label{tobeweakened}
\end{align}
In particular, for any $t>0$
the solution map $\Schwartz\owns (u(0), u_t(0)) \mapsto (u(t), u_t(t))$
for the Cauchy problem \eqref{NLW} 
fails to be continuous at $0$ in the $H^s \times H^{s-1}$ topology.

For $s< s_c$, the same is true, except that
the final conclusion \eqref{tobeweakened} is weakened to 
\begin{equation}
\|u(t)\|_{H^s} \ge c(t) \ \text{ for some constant $c(t)>0$
independent of $\eps$.}
\end{equation}
For $s=s_c$, the same holds with $c(t)$ independent of $t$.
Thus for $s\le s_c$, for any $t>0$, the solution operator 
fails to be continuous in $H^s$ at $0$.
\end{theorem}
\noindent Note that for $s<s_c$, Theorem~\ref{thm:bignorm-wave} gives
a stronger conclusion for spatial dimension $d=1$ if $s>0$ or $s<-\tfrac12$.

Theorem~\ref{thm:bignorm-1d}
implies the same conclusion in all higher dimensions.
\begin{corollary} \label{1dcorollary}
Suppose either that $p$ is an odd integer, or that
$p > \lfloor s \rfloor + 1$.
Then \eqref{NLW} is illposed, with rapid norm inflation in the sense
of Theorem~\ref{thm:bignorm-1d},
in $H^s(\reals^d)$ for all $d>1$ and $s<\max(0,\tfrac12-\tfrac1p)$.
\end{corollary}

\begin{figure}[htbp] \centering
 \ \psfig{figure=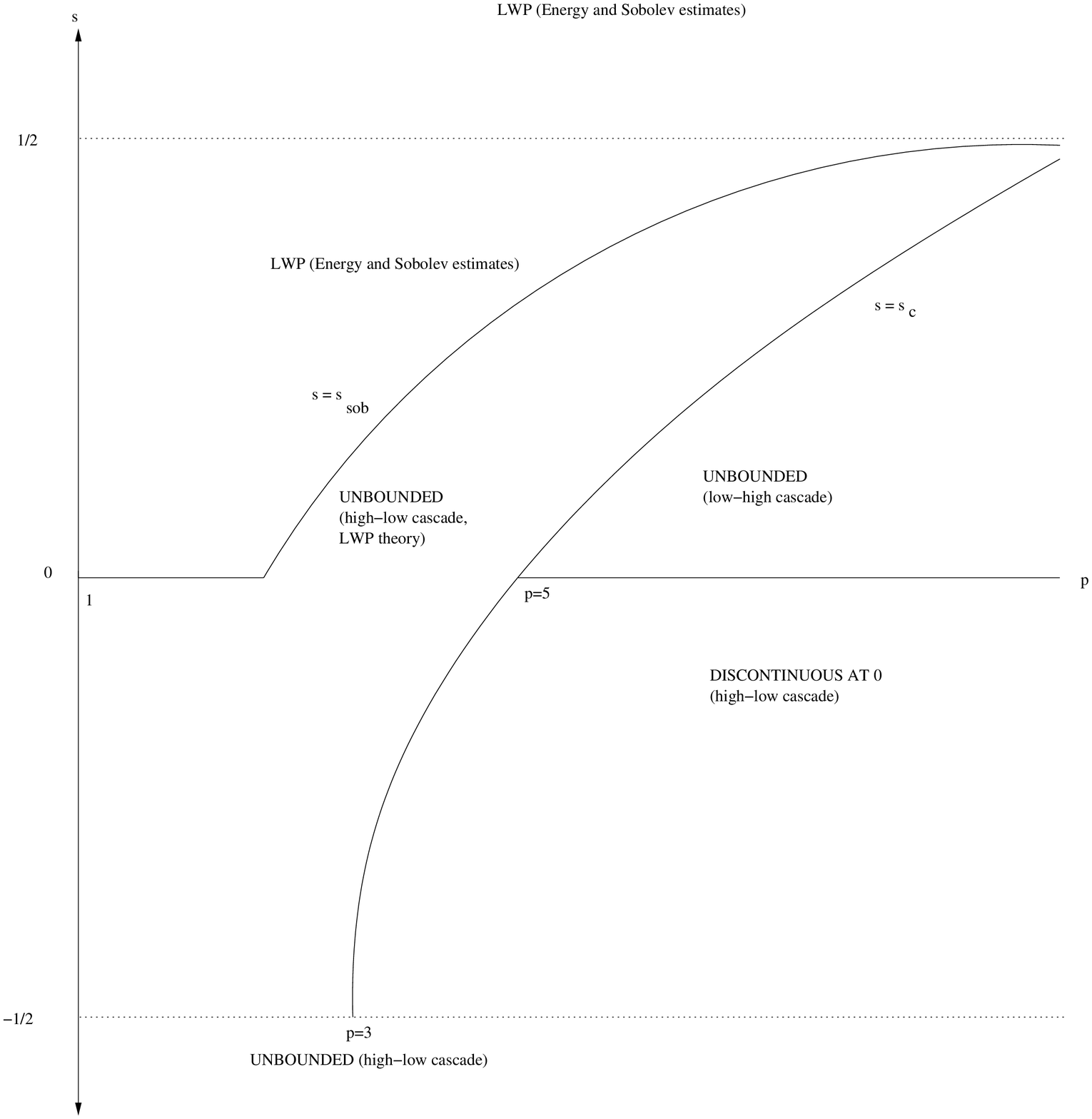,height=4in,width=5in}
 \caption{A schematic depiction of the NLW results for $d = 1$.  Here it is the Sobolev threshold
$s_{\rm sob}$ which dominates, rather than the scaling threshold $s_c$ or the Lorentz-invariance 
threshold $s_{\rm conf}$ (which always lies between $s_{\rm sob}$ and $s_c$).  
For $s \ge s_{\rm sob}$ well-posedness holds, by energy methods and Sobolev embedding. 
But for $s < s_{\rm sob}$ a high-to-low frequency cascade
causes rapid growth of $H^s$ norms. } 
 \end{figure}

We prove Theorem~\ref{thm:bignorm-1d}
in Section \ref{1d:sec}; the idea is to take initial data equal to an 
$H^{s_{\rm sob}}$ normalized approximate delta function and compute iterates, using the $H^{s_{\rm sob}}$ local well-posedness theory to control the convergence of the iterates.  There is a high-to-low cascade which moves the 
$H^{s_{\rm sob}}$ norm down from high frequencies to low frequencies, causing $H^s$ norm weak blowup.

\subsection{Further remarks}

In certain cases one obtains a stronger version of Theorems~\ref{thm:bignorm-1d}
and \ref{thm:bignorm-wave}:
(gNLW) evolutions explode instantaneously in all supercritical norms.
\begin{theorem}
  \label{thm:biggerbang}
Consider any parameters $s,p,d$ for which all of the following hold:
\newline
(i) $H^s$ norm inflation holds for (gNLW) in the sense of Theorems
\ref{thm:bignorm-wave} and \ref{thm:bignorm-1d}. 
\newline
(ii) For any smooth, compactly supported initial datum, there exists
a solution of (gNLW) in $C(\reals,H^s\cap C^\infty(\reals^d))$.
\newline
(iii) These solutions obey finite speed of propagation.
\newline
Then
there exist an initial datum $u(0) \in C^\infty \cap H^s$ and a
    corresponding solution $u(t) \in C^\infty$ for $t \in [0,1]$ such that
\[
\|u(t)\|_{H^s}=\infty \ \text{for all $t>0$.}
\]
\end{theorem}
By finite speed of propagation we mean that $u(t,x)=0$ if the initial data
vanish in a neighborhood of the closed ball of radius $t$ centered at $x$.
The hypotheses of Theorem~\ref{thm:biggerbang} 
are all satisfied for instance for defocusing equations for odd integers $p$
in spatial dimension $1$, and for the cubic defocusing equation in dimension
$d=2$, which is globally wellposed in $H^1$.

The proof of Theorem \ref{thm:biggerbang} is an adaptation of a
construction which is well-known (see e.g.\ \cite{sogge:wave}) in the
focusing setting. The inflationary solutions obtained in Theorems
\ref{thm:bignorm-wave} and \ref{thm:bignorm-1d} depend on a small parameter
$\epsilon$. We add up appropriate translates of solutions
$u_j$ depending on $\epsilon_j$, where $\epsilon_j \rightarrow 0$. 
By finite speed of propagation, the translates may be
arranged so that there is no interaction among the $u_j$. An appropriate choice of
$\{\epsilon_j\}$ produces an infinite norm when the contributions of the
pieces are summed. We omit the familiar details.

If the hypotheses of finite propagation speed and global existence for large smooth
data are dropped then one concludes that uniqueness, finite speed of
propagation, and existence in $H^s$ cannot simultaneously hold true.

All of our ill-posedness results for NLW continue to hold if a mass term $mu$ is
added to the equation, thus transforming it into a nonlinear Klein-Gordon equation; the arguments for \eqref{NLW} apply with small modifications (in particular, one has to compensate for the fact that the scaling $u \mapsto u^\lambda$ will also affect the mass $m$).
The point is that all of our examples are ``high-frequency'', and consequently
the mass term plays no significant role after rescaling.  
An alternative way of saying this is that our examples always have large $L^\infty$ 
norm, and so the nonlinear term $|u|^{p-1} u$ dominates the mass term.

We are grateful to Mark Keel for pointing out the connection with the work of Kuksin \cite{kuksin}.

\section{NLS: Small dispersion analysis}\label{section-zero}

The common element in all our arguments for \eqref{gNLSp} is
a quantitative analysis of the NLS equation 
\eqref{smalldispersion} with dispersion coefficient $\nu$
in the {\em small dispersion} regime $\nu \to 0$.  
Formally, as $\nu \to 0$ this equation approaches the ODE
\begin{equation} \label{smalldispersion-ode}
\left\{
\begin{aligned}
-i\phi_s(s,y) &= \omega |\phi|^{p-1} \phi(s,y)
\\
 \phi(0,y)&=\phi_0(y)
\end{aligned}
\right.
\end{equation}
which has the explicit solution $\phi = \phi^{(0)}$ defined by
\be{phi-ode}
\phi^{(0)}(s,y) := \phi_0(y) \exp(i \omega s |\phi_0(y)|^{p-1}).
\end{equation}

Whereas the standard well-posedness theory treats the the nonlinear Schr\"odinger
equation as a perturbation of the linear Schr\"odinger equation,
we take the opposite point of view and regard the nonlinearity 
$\omega|\phi|^{p-1} \phi$ as the main term, and the dispersive term $\nu^2 \Delta_y \phi$
as the perturbation.
When $\nu$ is nonzero but small, $\phi$ may be expected 
to stay close to $\phi^{(0)}$, at least for short times. 
A quantitative statement is as follows.
\begin{lemma}  \label{energylemma}
Let $d\ge 1$, $p \geq 1$, and let $k > d/2$ be an integer.  
If $p$ is not an odd integer, then we assume also the additional regularity condition 
$ p \geq k+1$.
Let $\phi_0$ be a Schwartz function.  Then there exists  
$C, c > 0$ depending on all the above parameters, such that if $0 < \nu \leq c$ 
is a sufficiently small real number, then for $T = c |\log \nu|^c$
there exists a solution $\phi(s,y)\in C^1([-T,T],H^{k,k})$ 
of \eqref{smalldispersion} satisfying
\begin{equation}\label{sd}
\|\phi(s) - \phi^{(0)}(s)\|_{H^{k,k}(\reals^d)}
\le C \nu \ \text{ for all } |s|\le  c |\log \nu|^c,
\end{equation}
where $H^{k,k}$ denotes the weighted Sobolev space
$$ \| \phi \|_{H^{k,k}} := \sum_{j = 0}^k
\| (1 + |x|)^{k-j} \partial_x^j \phi \|_{L^2}.$$
\end{lemma}

The proof is a straightforward application of the energy method; 
for the sake of completeness we sketch a proof.

\begin{proof}  We define the function $F: \complex \to \complex$ by 
\be{F-def}
F(z) := \omega |z|^{p-1} z,
\end{equation}
thus 
$$ -i \partial_s \phi^{(0)} = F(\phi^{(0)})$$
and the equation to be solved is
$$ -i \partial_s \phi + \nu^2 \Delta_y \phi = F(\phi).$$
Thus, with the Ansatz
\be{ansatz}
\phi = \phi^{(0)} + w,
\end{equation}
$w$ is a solution of the Cauchy problem 
\begin{equation*}
\left\{
\begin{aligned}
-i \partial_s w + \nu^2 \Delta_y w 
&= -\nu^2 \Delta_y (\phi^{(0)}) + F(\phi^{(0)}+w) - F(\phi^{(0)})\\
w(0,y) &= 0.
\end{aligned}
\right.
\end{equation*}
Our assumptions on $p$ guarantee that $F$ is a $C^k$ function,
whose first $k$ derivatives are all Lipschitz.  Since $k > d/2$,
the weighted Sobolev space $H^{k,k}(\reals^d)$ controls $L^\infty$,
and standard energy method arguments\footnote{The space
  $H^{k,k}$ is preserved by the free Schr\"odinger flow locally in time (see
  \cite{cct1} for a proof), so there is no difficulty in applying the
  energy method.} \cite{CazenaveBook} 
shows that a unique $H^{k,k}$ solution $w$ to the above Cauchy problem exists 
locally in time, and can be extended as long as $\| w(t) \|_{H^{k,k}}$ 
stays bounded.  Thus it will suffice to show that
\be{ws-hk-bound}
\sup_{|s| \leq T} \| w(s) \|_{H^{k,k}(\reals^d)} \leq C \nu
\end{equation}
where $0 \leq T \leq c |\log \nu |^c$ and we have the \emph{a priori} assumption that $w$ stays in $H^{k,k}$ for times $|s| \leq T$.

Without loss of generality we may take $s \geq 0$.  The energy inequality gives
$$ \partial_s \| w(s) \|_{H^{k,k}} 
\leq C \big\| -\nu^2 \Delta_y (\phi^{(0)})(s) 
+ F(\phi^{(0)}+w)(s) - F(\phi^{(0)})(s) \big\|_{H^{k,k}}
+ C\|w(s)\|_{H^{k,k}}.
$$
Since $\phi_0$ is Schwartz and $F$ is in $C^k$, we see from \eqref{phi-ode} that
\be{phiode-ck}
\| \phi^{(0)}(s) \|_{H^{k,k}} + \| \phi^{(0)}(s) \|_{C^k} \leq C(1+|s|)^k
\end{equation}
and then
$$ \| \Delta_y (\phi^{(0)})(s) \|_{H^{k,k}} \leq C (1 + |s|)^k$$
where $C$ depends on $\phi_0$, $k$, $p$, $d$.  

We now consider pointwise bounds for the quantity
$$ F(\phi^{(0)}+w)(s) - F(\phi^{(0)})(s).$$ From the 
mean value theorem and direct computation follow bounds
$$ |F^{(j)}(z)| \leq C (1 + |z|)^{p-j}$$
and
$$ |F^{(j)}(z+z') - F^{(j)}(z)| \leq C |z'|(1 + |z| + |z'|)^{p-j-1}$$
for $j = 0, 1, \ldots, k$.  Multiple applications of the chain rule, 
the Leibnitz rule, and \eqref{phiode-ck} thus imply the pointwise inequality
$$ |F^{(j)}(\phi^{(0)}+w)(s,y) - F^{(j)}(\phi^{(0)})(s,y)|
\leq C (1+|s|)^C \sum_{\alpha_1, \ldots, \alpha_r} 
|w^{(\alpha_1)}(s,y)| \ldots |w^{(\alpha_r)}(s,y)| (1 + |w(s,y)|^{p-r})$$
for each $0\le j\le k$, where $\alpha_1, \ldots, \alpha_r$ range over all finite collections of non-negative integers with $1 \leq r \leq k+1$ and $\alpha_1 + \ldots + \alpha_r \leq k$.
H\"older's and Sobolev's inequalities thus yield the bounds
$$ \| F(\phi^{(0)}+w)(s) - F(\phi^{(0)})(s) \|_{H^{k,k}}
\leq C (1+|s|)^C \|w(s) \|_{H^{k,k}} (1 + \| w(s) \|_{H^{k,k}})^{p-1}.$$
All together these estimates give the differential inequality
$$ \partial_s \| w(s) \|_{H^{k,k}} \leq C (1 + |s|)^C \nu^2 + C (1 + |s|)^C \| w(s) \|_{H^{k,k}} + C(1+ |s|)^C \| w(s) \|_{H^{k,k}}^p.$$
Under the {\it a priori} 
assumption that $w(s)$ is bounded in $H^{k,k}$, e.g.\ $\|w(s) \|_{H^{k,k}} \leq 1$, 
this becomes 
$$ \partial_s \| w(s) \|_{H^{k,k}} \leq C (1 + |s|)^C \nu^2 + C (1 + |s|)^C \| w(s) \|_{H^{k,k}}$$
which by Gronwall's inequality and the initial condition $w(0) = 0$ implies that
$$ \| w(s) \|_{H^{k,k}} \leq C \nu^2 \exp(C(1 + |s|)^C).$$
Thus if $|s| \leq c |\log \nu|^c$ for suitably chosen $c$ and $\nu$ is sufficiently small, we obtain 
\eqref{ws-hk-bound}, and furthermore we can recover the {\it a priori} 
assumption $\| w(s) \|_{H^{k,k}} \leq 1$, which can then be removed by the usual 
continuity argument.  This concludes the proof of the Lemma.
\end{proof}

\begin{remark} Most of our applications require merely $H^k$ control, rather than
$H^{k,k}$; the decay bounds are however useful for obtaining uniform bounds 
on Fourier coefficients.  From \eqref{phiode-ck} and Lemma \ref{energylemma} 
we obtain in particular the bounds 
\be{phi-bounds}
\| \phi(s) \|_{H^k} \leq C |\log |\nu||^C \hbox{ for all } |s| \leq c |\log |\nu||^c.
\end{equation}
\end{remark}

\begin{remark}
Lemma~\ref{energylemma}
is stated only for real $\nu$, but it is possible to add a small negative 
(respectively positive) imaginary part to $\nu$, thus introducing some dissipation into 
\eqref{smalldispersion}, and still obtain the same result for positive 
(respectively negative) times $s$.  We shall not pursue this matter here.
\end{remark}

\begin{remark}
It might be interesting to analyze the solution for a much longer
timespan than $T=c |\log \nu|^c$, but any $T$ which tends to infinity
as $\nu \to 0$ suffices for our applications. The control obtained here
should be contrasted with more global control obtained in \cite{cct1} 
for the 1D cubic NLS equation, which used a different approximating ODE 
and also required some smallness conditions on the data. 
\end{remark}

\section{NLS decoherence} \label{section:negatives}

We now use Lemma \ref{energylemma} to prove Theorem \ref{thm:NLS} in the case $s < 0$.  We wish to demonstrate failure of uniform continuity of the solution map.  Our argument 
is similar to those in \cite{cct1}, except that small dispersion solutions are used
in instead of modified scattering solutions.

Fix $s<0$.  Let $w(x)$ be an arbitrary nonzero Schwartz function. 
Let $a\in[\tfrac12,1]$ and
$\nu\in(0,1]$ be parameters; we will construct a family of
solutions of \eqref{gNLSp} depending on $a,\nu$, and will analyze them
quantitatively as $\nu\searrow 0$.

Let $k>d/2$ be an integer, and (if $p$ is not an odd integer) let $k+1 \leq p$. From Lemma~\ref{energylemma}, it follows that for $\nu \leq c$ there exists a solution $\phi(s,y) = \phi^{(a,\nu)}(s,y)$ to the equation \eqref{smalldispersion} with initial datum 
\be{phi-init}
\phi^{(a,\nu)}(0,y) := aw(y).
\end{equation}
  Note that the parameter $a$ varies within a compact set, and so the constants 
$C$, $c$ in Lemma~\ref{energylemma} may be taken to be independent of $a$.

Lemma \ref{energylemma}
gives the estimates
\be{phi-close}
\| \phi^{(a,\nu)}(s) - \phi^{(a,0)}(s) \|_{H^k(\reals^d)} \leq C \nu 
\hbox{ for all } |s| \leq c |\log \nu|^c
\end{equation}
where
\be{phia0-def}
\phi^{(a,0)}(s,y) := a w(y) \exp(i \omega a^{p-1} s |w(y)|^{p-1}).
\end{equation}
Let $0 < \delta \ll 1$, and suppose that $a, a' \in [1/2,1]$ are such that 
$|a-a'| \sim \delta$.  Then $\phi^{(a,0)}(s)$ and $\phi^{(a',0)}(s)$ 
differ by $O(\delta)$ at time $s=0$. 
On the other hand, if $s$ is sufficiently large that 
$\big||a|^{p-1}s-|a'|^{p-1}s\big|\sim \delta s$ has order of magnitude $1$, 
and if $\nu$ is chosen small enough to ensure that 
$1 \ll |s| \leq c |\log \nu|^c$, 
these functions become well separated due to 
phase decoherence, as is apparent from  \eqref{phia0-def}.  
The idea is to rescale this solution so that the decoherence occurs much sooner, while also using the Galilean invariance to keep the $H^s$ norm under control for $s < 0$.  For this argument the dominant parameter will be the Galilean invariance parameter $v$ (which can be used to shrink the $H^s$ norm for $s < 0$ by sending $|v|$ to infinity); the other parameters play only a supporting role, serving mainly to make the decoherence occur at arbitrarily small times.

We turn to the details. From \eqref{u-phi}, \eqref{nls-scaling}, \eqref{gal} we obtain a four-parameter family of solutions $u = u^{(a,\nu,\lambda,v)}$ to the NLS equation \eqref{gNLSp}, where $a \in [\frac{1}{2},1]$, $0 < \nu \ll 1$, $0 <\lambda \ll 1$, $v \in \R^d$, and
\be{u-def}
u^{(a,\nu,\lambda,v)}(t,x) :=
\lambda^{-2/(p-1)}
e^{-iv \cdot x/2} e^{i |v|^2 t/4}
\phi^{(a,\nu)}(\lambda^{-2} t, \lambda\rp\nu (x-vt)).
\end{equation}

Informally, for times $|t| \ll \lambda^2 |\log \nu|^c$, the function $u(t)$ has magnitude $\sim \lambda^{-2/(p-1)}$ and oscillates with frequency $\sim |v|$, with respect to the
spatial variable, on a ball of radius $\sim \lambda/\nu$.  
Thus one expects the $H^s$ norm to be of the order of $\lambda^{-2/(p-1)} |v|^s 
(\lambda/\nu)^{d/2}$.  The following lemma makes that intuition precise.

\begin{lemma}\label{u-estimates}  
Let $0\ne w\in\Schwartz$ be fixed.
Let $s<0$, and suppose that $k, \phi,\phi^{(0)}$ 
satisfy \eqref{sd}.
Suppose that $|v| \geq 1$, $a,a' \in [1/2,2]$, and
$0 < \nu \leq \lambda \ll 1$.  Then 
\be{a-bound}
\| u^{(a,\nu,\lambda,v)}(0) \|_{H^s} \leq C \lambda^{-2/(p-1)} |v|^s (\lambda/\nu)^{d/2}
\end{equation}
and
\be{aa-bound}
\| u^{(a,\nu,\lambda,v)}(0) - u^{(a',\nu,\lambda,v)}(0)\|_{H^s} \leq C \lambda^{-2/(p-1)} |v|^s (\lambda/\nu)^{d/2} |a-a'|.
\end{equation}
Likewise
\begin{multline} \label{t-diff}
 \| u^{(a,\nu,\lambda,v)}(t) - u^{(a',\nu,\lambda,v)}(t) \|_{H^s} 
\\
\geq
c |v|^s \lambda^{-2/(p-1)} (\lambda/\nu)^{d/2}
\Big(\| \phi^{(a,\nu)}(t/\lambda^2) - \phi^{(a',\nu)}(t/\lambda^2) \|_2 
- C |\log \nu|^C (\lambda/\nu)^{-k}|v|^{-s-k}\Big)
\end{multline}
whenever $|t| < c |\log \nu|^c\lambda^2$.
\end{lemma}

\begin{proof}
Defining the Fourier transform by
$$ \hat u(\xi) :=  \int_{\R^d} e^{-i x \cdot \xi} u(x)\ dx,$$
we have
$$ [u^{ (a,\nu,\lambda,v) }(0)]\widehat{\phantom{A}}(\xi) =
\lambda^{-2/(p-1)} ({\lambda}/{\nu})^d 
[\phi^{(a,\nu)}(0)]\widehat{\phantom{A}}(\frac{\lambda}{\nu}(\xi + v/2)).$$
By Plancherel, the change of variables 
$\eta = {\lambda}{\nu\rp}(\xi + v/2)$,  and the definition $\phi^{(a,\nu)}(0)\equiv aw$,
$$ \| u^{(a,\nu,\lambda,v)}(0) \|_{H^s}^2 \sim
\lambda^{-4/(p-1)} ({\lambda}/{\nu})^{d}
\int (1 +  |\frac{\nu}{\lambda} \xi-v/2|)^{2s} 
|a\widehat{w}(\xi)|^2\ d\xi.$$
Since $|a|\sim 1$, 
\begin{align*}
 \| u^{(a,\nu,\lambda,v)}(0) \|_{H^s}^2
&\le C \lambda^{-4/(p-1)} (\lambda/\nu)^d 
|v|^{2s} \int_{|\xi|\le \lambda\nu\rp|v|/4} |\widehat{w}(\xi)|^2\,d\xi
\\
&\qquad\qquad + C \lambda^{-4/(p-1)} (\lambda/\nu)^d  
\int_{|\xi|\ge \lambda\nu\rp|v|/4} |\widehat{w}(\xi)|^2\,d\xi
\\
&\le C \lambda^{-4/(p-1)} (\lambda/\nu)^d |v|^{2s}
+ C_N \lambda^{-4/(p-1)} (\lambda/\nu)^d (\lambda\nu\rp|v|)^{-N}
\end{align*}
for all $N<\infty$, since $w\in\Schwartz$.
Since $\lambda\ge\nu$, this is 
$\le C\lambda^{-4/(p-1)} (\lambda/\nu)^d |v|^{2s}$,
as desired.
Since $u^{(a,\nu,\lambda,v)}(0)- u^{(a',\nu,\lambda,v)}(0) \equiv (a-a')
 u^{(1,\nu,\lambda,v)}(0)$, we likewise have
$$ \| u^{(a,\nu,\lambda,v)}(0)- u^{(a',\nu,\lambda,v)}(0) \|_{H^s} 
\leq C
\lambda^{-2/(p-1)} ({\lambda}/{\nu})^{d/2} |v|^{s} |a-a'|\|w \|_{H^k}.$$
The claims \eqref{a-bound}, \eqref{aa-bound} follow.

Now consider the same computations at a nonzero time $t$.  
The phase factor $e^{i |v|^2 t/4}$ in \eqref{u-def} is harmless, 
as is the translation $x-vt$, whence 
\begin{multline*}
\| u^{(a,\nu,\lambda,v)}(t) - u^{(a',\nu,\lambda,v)}(t) \|_{H^s}^2 
\\
\sim
\lambda^{-4/(p-1)} (\lambda/\nu)^{d} 
\int (1 +  |\nu\lambda\rp \xi-v/2|)^{2s} 
|[\phi^{(a,\nu)}(t/\lambda^2)]\widehat{\phantom{A}}(\xi)
- [\phi^{(a',\nu)} (t/\lambda^2)]\widehat{\phantom{A}}(\xi)|^2\ d\xi.
\end{multline*}
In particular there is the lower bound
\begin{multline*}
\| u^{(a,\nu,\lambda,v)}(t) - u^{(a',\nu,\lambda,v)}(t) \|_{H^s}^2 
\\
\geq
c \lambda^{-4/(p-1)} (\lambda/\nu)^{d}|v|^{2s}
\int_{|\xi| \leq |v| \lambda / 4\nu}
|[\phi^{(a,\nu)}(t/\lambda^2)]\widehat{\phantom{A}}(\xi)
- [\phi^{(a',\nu)}(t/\lambda^2)]\widehat{\phantom{A}}(\xi)|^2\ d\xi.
\end{multline*}
On the other hand,
$$ \int_{|\xi| \geq |v| \lambda / 4\nu}
\big(|[\phi^{(a,\nu)}(t/\lambda^2)]\widehat{\phantom{A}}(\xi)|^2
+ |[\phi^{(a',\nu)}(t/\lambda^2)]\widehat{\phantom{A}}(\xi)|^2\big)\,d\xi
\leq C |\log \nu|^C (\lambda/\nu)^{-2k}|v|^{-2k}.$$
This follows from \eqref{phi-bounds} and Plancherel, since 
$\lambda/\nu$ and $|v|$ are $ \geq 1$, and $s\le 0$.
\end{proof}

To exploit Lemma~\ref{u-estimates} to prove Theorem~\ref{thm:NLS} we set
$$ \lambda := \nu^\sigma$$
for some small $\sigma = \sigma(p,s,d) > 0$ to be chosen later.  We then choose any
vector $v$ whose magnitude $|v|$ satisfies
$$ \lambda^{-2/(p-1)} |v|^s (\lambda/\nu)^{d/2} = \eps,$$
or equivalently
$$ |v| = \nu^{\frac{1}{s}(\frac{d(1-\sigma)}{2} + \frac{2\sigma}{p-1})} \eps^{1/s},$$
where $0 < \eps \ll 1$ is a parameter to be chosen later.
Since $s < 0$, we see that the power to which $\nu$ is raised
is negative if $\sigma$ is chosen to be sufficiently  small.  
In particular, $|v|$ grows faster than any power of $|\log \nu|$ as $\nu \to 0$;
this fact will be used repeatedly in the sequel.

Let $a$, $a'$ be arbitrary distinct elements of $[1/2,2]$.  
From Lemma~\ref{u-estimates} it follows that
$$
\| u^{(a,\nu,\lambda,v)}(0) \|_{H^s}\ +\  \| u^{(a',\nu,\lambda,v)}(0) \|_{H^s} \leq C \eps
$$
and
$$
\| u^{(a,\nu,\lambda,v)}(0) - u^{(a',\nu,\lambda,v)}(0)\|_{H^s} \leq C \eps |a-a'|.
$$
On the other hand, from an inspection of \eqref{phia0-def} we see that there exists a time $T = T(a,a') > 0$ such that
$$ \| \phi^{(a,0)}(T) - \phi^{(a',0)}(T)\|_{L^2} \geq c,$$
where $c>0$ is independent of $a,a'$.
Fix this $T$.  \eqref{phi-close} and the triangle inequality give
$$ \| \phi^{(a,\nu)}(T) - \phi^{(a',\nu)}(T)\|_{L^2} \geq c$$
if $\nu$ is sufficiently small and $T\le c|\log \nu |^c$.  
From \eqref{t-diff} we thus conclude
$$
\| u^{(a,\nu,\lambda,v)}(\lambda^2 T) - u^{(a',\nu,\lambda,v)}(\lambda^2 T)\|_{H^s} 
\geq c \eps
-C(\lambda/\nu)^{-k}|v|^{-s-k}|\log\nu|^C,
$$
provided that $T\le c|\log \nu |^c$,
again if $\nu$ is sufficiently small (since $|v|$ grows like a negative power of $\nu$).

If $\sigma$ is chosen to be sufficiently small then
$(\lambda/\nu)^{-k}|v|^{-s-k}|\log\nu|^C\to 0$ as $\nu\to 0$.
Indeed, 
\[
(\lambda/\nu)^{-k}|v|^{-s-k}
= \nu^{k-(s+k)d/2s + O(\sigma)},
\]
and since $s<0$ and $k>d/2$, the exponent is
\[
k(1-\frac{d}{2s})-\frac{d}2 >
\frac{d}2(1-\frac{d}{2s})-\frac{d}2 
= -\frac{d^2}{4s}>0.
\]
Therefore
$$
\| u^{(a,\nu,\lambda,v)}(\lambda^2 T) - u^{(a',\nu,\lambda,v)}(\lambda^2 T)\|_{H^s} 
\geq c \eps
$$
for all sufficiently small $\nu$.

Letting $\nu\to 0$ (and hence $\lambda\to 0$) forces
the time $\lambda^2 T$ to approach zero as well.
Thus we have constructed two solutions whose $H^s$ norms are of size $O(\eps)$ at time zero, 
and which are distance $O(\eps |a-a'|)$ apart in $H^s$ norm norm at time zero, 
but at some arbitrarily short time later are separated in $H^s$ norm
by at least $c\eps$.  Since $|a-a'|$ can be made arbitrarily small, 
this shows that the solution map\footnote{Note that the solutions $u^{(a,\nu,\lambda,v)}$ 
here have sufficient regularity to be covered by the uniqueness theory of \cite{cwI}.} 
cannot be uniformly continuous, as claimed.
\qed

\begin{remark}
In the case when $p$ is an odd integer, the restriction $\lambda\ge\nu$
in Lemma~\ref{u-estimates} is inessential. 
In that case, Lemma~\ref{energylemma} yields control of
solutions in $H^k$ for all finite $k$.
A more careful repetition of the proof of Lemma~\ref{u-estimates} then reveals
that its conclusions remain valid, provided
merely that there exist $\rho>0$ such that
$\lambda\nu\rp \ge c|v|^{-1+\rho}$, assuming always that $|v|\gg 1$ and $s<0$.
Equivalently, $\nu\lambda\rp \le C|v|^{\varrho}$ for some $\varrho<1$.

This in turn allows for greater flexibility in the choice of the parameter $\sigma$
in the application to ill-posedness.
Consider the fundamental case $d=1$, $p=3$ for simplicity.
Then the parameters $\lambda,\nu,v$ are related by
$\lambda = \nu^\sigma$ and $\nu = |v|^{s/(1+\sigma)}$,
with the constraint $\sigma>0$.
The initial data in our construction are parametrized by $|v|$, and have Fourier transforms 
supported where $|\xi+v/2|\le C\nu\lambda\rp = C|v|^\varrho$
where $\varrho = (-s)(\sigma-1)/(\sigma+1)$.
This is less than $-s$, 
and for $-1\le s<0$ it can be made arbitrarily close to $-s$
by choosing $\sigma$ arbitrarily large.
Thus rapid decoherence is seen for initial data whose Fourier transforms
are supported essentially throughout intervals of comparatively large diameters 
$|v|^{|s|-\eps}$ centered at $v/2$, as $|v|\to\infty$. 
\end{remark}

\section{NLS energy transfer: The case $ 0< s < s_c$} 
\label{section:positives}

This section is devoted primarily to the proof 
of Theorem~\ref{thm:bignorm} in the case $0< s < s_c$.  Note that this theorem 
directly implies Theorem~\ref{thm:NLS} for the same range of values of $s$.

We again use the family $u^{(a,\nu,\lambda,v)}$ of solutions from the previous section.  
However, the parameter $v$ will not be used here, and we simply set $v := 0$; 
now the scaling parameter $\lambda$ will do most of the work.

Suppose that  $a \in [1/2,2]$.  From \eqref{u-def}, \eqref{phi-init} we have
$$ u^{(a,\nu,\lambda,0)}(0,x) = 
\lambda^{-2/(p-1)} a w({\nu x}/{\lambda}).$$
In particular, 
$$  [u^{(a,\nu,\lambda,0)}(0)]\widehat{\phantom{A}}(\xi) =
a\lambda^{-2/(p-1)} ({\lambda}/{\nu})^d \widehat{w}({\lambda \xi}/{\nu})$$
and hence
\begin{align*}
\|u^{(a,\nu,\lambda,0)}(0)\|_{H^s}^2 
&= a^2\lambda^{-4/(p-1)} (\lambda/\nu)^{2d}
\int |\widehat{w}(\lambda\nu\rp\xi)|^2 (1+|\xi|^2)^s\,d\xi
\\
&= a^2 \lambda^{-4/(p-1)} (\lambda/\nu)^{d}
\int |\widehat{w}(\eta)|^2 (1+|\nu\lambda\rp\eta|^2)^s\,d\eta.
\\
&\sim
\lambda^{-4/(p-1)} (\lambda/\nu)^{d-2s}
\int_{|\eta|\ge\lambda\nu\rp} |\widehat{w}(\eta)|^2|\eta|^{2s}\,d\eta
\\
& \qquad\qquad+
\lambda^{-4/(p-1)} (\lambda/\nu)^{d}
\int_{|\eta|\le\lambda\nu\rp} |\widehat{w}(\eta)|^2\,d\eta
\\
&=
\lambda^{-4/(p-1)} (\lambda/\nu)^{d-2s}
\int_{\reals^d} |\widehat{w}(\eta)|^2|\eta|^{2s}\,d\eta
\\
&\qquad\qquad- \lambda^{-4/(p-1)} (\lambda/\nu)^{d-2s}
\int_{|\eta|\le\lambda\nu\rp} |\widehat{w}(\eta)|^2\,
\big((\lambda/\nu)^{2s}-|\eta|^{2s}\big)\,d\eta.
\end{align*}
Thus for any $s>-d/2$,
\begin{equation}
\|u^{(a,\nu,\lambda,0)}(0)\|_{H^s}
=  c\lambda^{-2/(p-1)}({\lambda}/{\nu})^{d/2 - s}
\cdot(1 + O\big((\lambda\nu\rp)^{s+d/2}\big))
\end{equation}
where $c\ne 0$ provided that $w$ is not identically zero.
In particular, 
\begin{equation} \label{time0normal}
\|u^{(a,\nu,\lambda,0)}(0)\|_{H^s}
\le C\lambda^{-2/(p-1)}({\lambda}/{\nu})^{d/2 - s}
= C\lambda^{s_c-s}\nu^{s-d/2}
\end{equation}
provided that $s>-d/2$ and $\lambda\le\nu$.

For $s\le -d/2$, \eqref{time0normal} still holds, under the supplementary
hypothesis that 
\begin{equation}  \label{momentcondition}
\widehat{w}(\xi)
=O(|\xi|^\kappa) \ \text{as $\xi\to 0$, for some $\kappa> -s-d/2$.}
\end{equation}
For then, if $\lambda\le\nu$,
$\int_{\reals^d} |\widehat{w}(\eta)|^2|\eta|^{2s}\,d\eta<\infty$
and
\begin{equation*}
\int_{|\eta|\le\lambda\nu\rp} |\widehat{w}(\eta)|^2\,
\big((\lambda/\nu)^{2s}-|\eta|^{2s}\big)\,d\eta
\le C(\lambda\nu\rp)^{d+2s+2\kappa}
\le C<\infty.
\end{equation*}
This variant will be used in \S\ref{section:verylows}.

Consider now any $s<s_c$.
Let $\eps\in(0,1]$ be fixed for the time being, and set
\be{lambda-fix}
\lambda^{s_c - s} \nu^{s-d/2} = \eps;
\end{equation}
equivalently,
\[\lambda=\nu^\sigma\ \text{ where } \sigma = \frac{\frac{d}2-s}{s_c-s}>1.\]
Then $0 \leq \lambda \leq \nu$ as $\nu \to 0$ (since $s < s_c < d/2$), and 
$$  \| u^{(a,\nu,\lambda,0)}(0) \|_{H^s} \leq C \eps.$$

Specialize further to the case $0<s<s_c$. 
We investigate the behavior of $u^{(a,\nu,\lambda,0)}(t)$ for $t>0$,
starting with the analysis of $\phi^{(a,0)}(t,x)$ for $t \gg 1$.  
\eqref{phia0-def} gives
$$ \partial_x^j \phi^{(a,0)}(t,x) := a w(x) t^j (i \omega a^{p-1} \nabla_x |w(x)|^{p-1})^j \exp(i \omega a^{p-1} t |w(x)|^{p-1}) + O(t^{j-1})$$
for all integers $j\ge 0$ if $p$ is an odd integer, and all $j$ satisfying 
$0\le j\le p-1$ otherwise.
Since the support of $\phi^{(a,0)}(t,x)$ is independent of $t$,
this means that for all large enough $t$, depending on $j$, 
$$ \| \phi^{(a,0)}(t) \|_{H^j} \sim t^j$$
for the same values of $j$,
provided of course that $w$ does not vanish identically.
In particular, since the Sobolev norms $H^s$ are log-convex, 
$$ \| \phi^{(a,0)}(t) \|_{H^s} \sim t^s$$
whenever $s\ge 0$ is no larger than the greatest integer $\le p-1$.
If $p$ is not an odd integer,
then $s < s_c < d/2 < k$, and $p-1\ge k$; thus this conclusion holds for all $s$
under consideration in Theorem~\ref{thm:NLS}.
If $\nu \ll 1$ and $1 \ll t \leq c |\log \nu|^c$, \eqref{sd} thus implies that
\be{phi-growth}
\| \phi^{(a,\nu)}(t) \|_{H^s} \sim t^s.
\end{equation}

This estimate indicates that as time progresses, 
the function $\phi^{(a,\nu)(t)}$ transfers its energy to increasingly higher frequencies. 
 
We now exploit the supercriticality of $s$ via the scaling parameter $\lambda$ to 
create arbitrarily large $H^s$ norms at arbitrarily small times.
Applying \eqref{u-def} as before, we have
$$ [ u^{(a,\nu,\lambda,0)}(\lambda^2 t)]\widehat{\phantom{A}}(\xi) =
\lambda^{-2/(p-1)} ({\lambda}/{\nu})^d [\phi^{(a,\nu)}(t)]\widehat{\phantom{A}}
({\lambda \xi}/{\nu})$$
and hence, by the change of variables $\eta := \lambda \xi/\nu$,
$$\| u^{(a,\nu,\lambda,0)}(\lambda^2 t) \|_{H^s}^2
\geq c \lambda^{-4/(p-1)} ({\lambda}/{\nu})^{d}
\int (1 + |{\nu\eta}/{\lambda}|)^{2s}
|[\phi^{(a,\nu)}(t)]\widehat{\phantom{A}}(\eta)|^2\ d\eta,$$
since $\lambda/\nu\le 1$.

Now
\begin{align*}
\int (1 + |{\nu\eta}/{\lambda}|)^{2s}
|[ \phi^{(a,\nu)}(t) ]\widehat{\phantom{A}}(\eta)|^2\,d\eta
&\ge (\nu/\lambda)^{2s} \int_{|\eta|\ge 1} |\eta|^{2s}|
[ \phi^{(a,\nu)}(t) ]\widehat{\phantom{A}}(\eta)|^2\ d\eta
\\
&\ge (\nu/\lambda)^{2s}\Big( 
c\|\phi^{(a,\nu)}(t)\|_{H^s}^2
- C\|\phi^{(a,\nu)}(t)\|_{H^0}^2\big).
\end{align*}
From \eqref{phi-growth} it is apparent that $\|\phi^{(a,\nu)}(t)\|_{H^0}
\ll \|\phi^{(a,\nu)}(t)\|_{H^s}$ for $|t|\gg 1$.
Thus
$$ \| u^{(a,\nu,\lambda,0)}(\lambda^2 t) \|_{H^s}
\geq c \lambda^{-2/(p-1)} ({\lambda}/{\nu})^{d/2-s}
\| \phi^{(a,\nu)}(t) \|_{H^s} \geq c \eps t^s$$
by \eqref{lambda-fix}, \eqref{phi-growth}.  Theorem \ref{thm:bignorm} then follows, 
in the case $0 < s < s_c$,  by choosing $t$ large enough depending on $\eps$, 
and then $\nu$ (and hence $\lambda$) small enough depending on $\eps$, $t$.  
Theorem~\ref{thm:NLS}, for the same range of $s$, follows directly. 
\qed

\medskip
One case of Theorem \ref{thm:NLS} remains outstanding, namely when $s=0 < s_c$.  
The $H^0$ norm is conserved for Schwartz class initial data, 
there is certainly no norm growth result for $H^0$ like that established above for
$0<s<s_c$. Nevertheless the failure of the solution operator to be uniformly
continuous can be proved by an adaptation of the preceding construction.

We merely sketch the argument.  
Let $\eps>0$ be an arbitrary small number.
Consider two distinct numbers $a, a' \in [1/2,2]$. 
Let parameters $\lambda,\nu$ be related as in \eqref{lambda-fix}.
A direct calculation using \eqref{phia0-def} shows that
$$ \|  \phi^{(a,0)}(t) - \phi^{(a',0)}(t) \|_{L^2} \geq c>0$$
for some time $t$ satisfying
$|a-a'|\rp\leq t \leq c |\log \nu|^c$, 
provided of course that $\nu$ is chosen to be so small that $c|\log\nu|^c\ge |a-a'|\rp$.
By \eqref{sd}, this implies 
$$ \|  \phi^{(a,\nu)}(t) - \phi^{(a',\nu)}(t) \|_{L^2} \geq c$$
again if $\nu$ is sufficiently small.  Applying \eqref{u-def}, \eqref{lambda-fix} as before, we eventually obtain
$$ \|  u^{(a,\nu,\lambda,0)}(t \lambda^2) - u^{(a',\nu,\lambda,0)}(t \lambda^2) \|_{L^2} \geq c \eps.$$
However, by arguing as in the previous section 
(or by direct computation using \eqref{u-def}, \eqref{lambda-fix}) we have
$$ \|  u^{(a,\nu,\lambda,0)}(0) \|_{L^2} \leq C \eps$$
and
$$ \|  u^{(a,\nu,\lambda,0)}(0) - u^{(a',\nu,\lambda,0)}(0) \|_{L^2} \leq C \eps |a-a'|.$$
Since $|a-a'|$ can be arbitrarily small, this contradicts 
uniform continuity of the solution map as before. The 
proof of Theorem \ref{thm:NLS} is now complete.
\qed

\section{NLS energy transfer: The case $s < -d/2$} \label{section:verylows}

We now prove Theorem \ref{thm:bignorm} in the case $s \le  -d/2$.  The argument
is rather similar to that given above for $0<s<s_c$, except that it relies on a 
transfer of energy from high modes to low modes, rather than from low to high.

Let $w$ and $u^{(a,\nu,\lambda,0)}$ be as in the previous section. 
As shown in \S\ref{section:positives},
$\| u^{(a,\nu,\lambda,0)}(0) \|_{H^s} \leq C \eps$
if $0 < \lambda \leq \nu \ll 1$ are chosen according to \eqref{lambda-fix},
provided that $w$ is chosen so that $\widehat{w}(\xi)=O(|\xi|^\kappa)$
as $\xi\to 0$, for some $\kappa>-s-d/2$; we assume this moment condition 
henceforth.

Although $\widehat{w} (0)=0$, the function 
$\phi^{(a,0)}(s,x) = aw(x)\exp(i\omega a^{p-1}|w(x)|^{p-1}s)$
introduced in \eqref{phia0-def} 
need not share this property. 
Indeed, it is clear that $w$, $a$ can be chosen so that
$$ |\int \phi^{(a,0)}(1,y)\ dy| \geq c$$
for some constant $c>0$.
Equivalently, 
$$ |\big[{\phi^{(a,0)}(1)}\big]\widehat{\phantom{A}}(0)| \geq c.$$
Since $\phi^{(a,0)}(1)$ is rapidly decreasing, we thus see by continuity that
$$ |[\phi^{(a,0)}(1)]\widehat{\phantom{A}}(\xi)| \geq c \hbox{ for } |\xi| \leq c,$$
with our convention that the value of $0 < c \ll 1$ is allowed to change 
freely from line to line.
On the other hand, from \eqref{sd} and Cauchy-Schwarz 
(since $H^{k,k}$ controls $L^1$ when $k > n/2$) we have 
$$ |[\phi^{(a,\nu)}(1)]\widehat{\phantom{A}}(\xi) 
- [\phi^{(a,0)}(1)]\widehat{\phantom{A}}(\xi)| \leq C \nu$$
and thus
\be{phi-large}
|[\phi^{(a,\nu)}(1)]\widehat{\phantom{A}}(\xi)| \geq c \hbox{ for } |\xi| \leq c,
\end{equation}
uniformly for all sufficiently small $\nu>0$.
Meanwhile, the definition \eqref{u-def} of $u^{(a,\nu,\lambda,0)}$ specializes
at time $1$ to
$$
u^{(a,\nu,\lambda,0)}(\lambda^2,x) =
\lambda^{-2/(p-1)}
\phi^{(a,\nu)}(1, \frac{\nu x}{\lambda})$$
and hence

$$ \big[ u^{(a,\nu,\lambda,0)}(\lambda^2)\big]\widehat{\phantom{A}}(\xi) =
\lambda^{-2/(p-1)}
(\lambda{\nu\rp})^d
[\phi^{(a,\nu)}(1)]\widehat{\phantom{A}}({\lambda\nu\rp \xi}),$$
thus we have
$$ |\big[ u^{(a,\nu,\lambda,0)}(\lambda^2)\big]\widehat{\phantom{A}}(\xi)|
\geq c \lambda^{-2/(p-1)} ({\lambda}{\nu\rp})^d$$
when $|\xi| \leq c \nu/\lambda$.  In particular, 
\begin{multline*} \| u^{(a,\nu,\lambda,0)}(\lambda^2) \|_{H^s}
\geq c \lambda^{-2/(p-1)} ({\lambda}{\nu\rp})^d
= c\lambda^{s_c}\lambda^{-d/2} (\lambda\nu\rp)^d
\\
= c(\lambda\nu\rp)^{s+(d/2)} \cdot\big(\lambda^{s_c-s}\nu^{(d/2)-s)} \big)
=c(\lambda\nu\rp)^{s+(d/2)}\eps
\end{multline*}
just by considering the region $|\xi| \leq c$.  From \eqref{lambda-fix} we thus have
$$ \| u^{(a,\nu,\lambda,0)}(\lambda^2) \|_{H^s}
\geq c \eps (\lambda/\nu)^{d/2 + s}.$$
As $\nu \to 0$, \eqref{lambda-fix} ensures 
that $\lambda/\nu \to 0$ and $\lambda \to 0$.  Since $s < -d/2$, 
$(\lambda/\nu)^{d/2 + s} \to \infty$.  The claim then follows by setting $\nu$ sufficiently small depending on $\eps$.
\qed

The endpoint $s=-d/2$ requires a small modification. 
We work now on the entire region $|\xi| \leq c \nu/\lambda$, 
instead of simply discarding everything outside the ball $|\xi| \leq c$. 
This leads to a bound 
$$  \| u^{(a,\nu,\lambda,0)}(\lambda^2) \|_{H^{-n/2}}
\geq c \eps \log (\lambda/\nu),$$
and the argument then proceeds as before. 
The remaining details are left to the reader.

\section{NLW: Small dispersion analysis}\label{nlw-zero-sec}


We now analyze the small dispersion approximation for
the NLW equation \eqref{NLW}.  In analogy with the NLS case, consider 
\begin{equation} \label{smalldispersion-wave}
\left\{
\begin{aligned}
-\phi_{ss}(s,y) + \nu^2 \Delta_y \phi(s,y) &= \omega |\phi|^{p-1} \phi(s,y)
\\
\phi(0,y)&=\phi_0(y)\\
\phi_t(0,y)&=0
\end{aligned}
\right.
\end{equation}
in the zero-dispersion limit\footnote{This limit is the opposite of the \emph{non-relativistic limit}, which sends $\nu$ to infinity!} $\nu \to 0$.  We have set the initial velocity equal to zero to avoid technicalities.
As in the Schr\"odinger case, for any solution $\phi$ of \eqref{smalldispersion-wave} ,
\begin{equation}  \label{wave-transform}
\lambda^{-2/(p-1)} \phi(\lambda\rp{t}, \lambda\rp{\nu x}),
\end{equation}
defines a solution of \eqref{NLW}.

Suppose $\phi_0$ is supported in the unit ball $|y| \leq 1$.  A standard finite speed of propagation argument (see e.g.\ \cite{sogge:wave}) then shows that $\phi(t)$ is supported in the ball $|y| \leq 1 + \nu |t|$ for all time $t$, provided that the solution stays regular 
(e.g.\ $\phi \in H^k$ for some $k > d/2$; actually one has finite speed of propagation 
for much rougher data, see e.g. \cite{sogge:wave}) up to that time $t$.
 
Formally setting $\nu=0$ in \eqref{smalldispersion-wave} gives the ODE
$$
\left\{
\begin{aligned}
-\phi_{ss}(s,y) &= |\phi|^{p-1} \phi(s,y)
\\
\phi(0,y)&=\phi_0(y)\\
\phi_t(0,y)&=0
\end{aligned}
\right.
$$
which has the explicit solution $\phi = \phi^{(0)}$, where
\be{phi0w-def}
\phi^{(0)}(s, y) := \F(|\phi_0(y)|^{(p-1)/2} s) \phi_0(y)
\end{equation}
and $\F: \R \to \R$ is the unique solution to the ODE
$$ -\F''(s) = |\F(s)|^{p-1} \F(s); \quad \F(0) = 1; \quad \F'(0) = 0.$$
This is the Hamiltonian flow on a two-dimensional phase space with Hamiltonian
$$ H := \frac{1}{2} |\F'(s)|^2 + \frac{1}{p+1} |\F(s)|^{p+1};$$
since the level surfaces of this Hamiltonian are closed curves, we see that $\F$ is a bounded non-constant periodic function.  Also, from the assumptions on $p$, $k$ it is easy to see that $\F$ is a $C^{k+2}$ function.  
$\F(s)$ may be viewed as a ``nonlinear'' variant of $\cos(s)$ 
(which is the $p=1$ version of this function).

The fact that $|\phi_0(y)|$ is raised to the power $(p-1)/2$ in
\eqref{phi0w-def} might cause some problems with smoothness, but we
will finesse this issue by requiring $\phi_0$ to be the square of a
smooth function.

Observe that $\phi^{(0)}$ is real-valued if $\phi_0$ is.
The analogue of Lemma \ref{energylemma} is

\begin{lemma}  \label{energylemma-wave}
Let $d\ge 1$, $p \geq 1$, and let $k > d/2$ be an integer.  If $p$ is not an odd integer, then we also assume the additional regularity assumption
$$ p \geq k+1.$$
Let $\phi_0$ be the square of a compactly supported $C^\infty$ function.  
Then there exist   $C, c > 0$ depending on all the above parameters, such that 
for each sufficiently small real number $0 < \nu \leq c$, there exists a solution $\phi(s,y)$ of 
\eqref{smalldispersion-wave} for all $|s|\le c |\log\nu|^c$ such that
\begin{equation}\label{sd-wave}
\|\phi(s) - \phi^{(0)}(s)\|_{H^{k}(\reals^d)}
+ \|\phi_s(s) - \phi^{(0)}_s(s)\|_{H^{k}(\reals^d)}
\le C |\nu| \ \text{ for all } |s|\le  c|\log\nu|^c.
\end{equation}
\end{lemma}

One could also place spatial weights in the $H^k$ norm as in Lemma \ref{energylemma}, but these are redundant thanks to finite speed of propagation.

\begin{proof}
We repeat the proof of Lemma \ref{energylemma}.  If we again make the Ansatz \eqref{ansatz}, then $w$ is a solution of the Cauchy problem
\begin{equation*}
\left\{
\begin{aligned}
-\partial_{ss} w + \nu^2 \Delta_y w &= -\nu^2 \Delta_y \phi^{(0)} + F(\phi^{(0)}+w) - F(\phi^{(0)})\\
w(0,y) &= 0\\
w_s(0,y) &= 0
\end{aligned}
\right.
\end{equation*}
where $F$ was given by \eqref{F-def}.

Again, the energy method (see \cite{sogge:wave}) shows that this problem has a local solution in $H^{k+1} \times H^{k}$, which exists as long as the $H^{k+1} \times H^k$ norm stays bounded.  (Note that the wave equation energy inequality
is smoothing of order 1, so one only needs $k$ degrees of regularity on the nonlinearity in order to place $w$ in $H^{k+1} \times H^k$.)

Define the $\nu$-energy of a wave $w$ by
$$ E_\nu(w(s)) := \int \frac{1}{2} |w_s(s,y)|^2 + \frac{\nu^2}{2} |\nabla_y w(s,y)|^2\ dy.$$
If we have $-\partial_{ss} w + \nu^2 \Delta_y w = F$, then the standard energy identity gives
$$ \partial_s E_\nu(w(s)) = - \int w_s(s,y) F(s,y)\ dy,$$
which by Cauchy-Schwarz thus gives
$$ |\partial_s E_\nu^{1/2}(w(s))| \leq C \| F(s) \|_2.$$
Similarly, if we define
$$ E_{\nu,k}(w(s)) := \sum_{j=0}^k E_\nu(\partial_y^j w(s))$$
then we have
\be{energy-ineq}
|\partial_s E_{\nu,k}^{1/2}(w(s))| \leq C \| F(s) \|_{H^k}
\end{equation}
since the operator $-\partial_{ss} + \nu^2 \Delta_y$ commutes with $\partial_y$.

At first glance it seems that the quantity $E_{\nu,y}$ does not give good control on the spatial derivatives of $w$, because of the factor of $\nu$.  
Fortunately, one can instead use the time derivatives to obtain regularity.  
Indeed, the fundamental theorem of calculus gives
\be{w-bound}
\| w(s) \|_{H^k} \leq \int_0^s \| w_s(s')\|_{H^k}\ ds'
\leq C \int_0^s E_{\nu,k}^{1/2}(w(s'))\ ds'
\leq C s e(s),
\end{equation}
where $e(s)$ is the non-decreasing function
$$ e(s) := \sup_{0 \leq s' \leq s} E_{\nu,k}^{1/2}(w(s')).$$
Also, since $\phi^{(0)}$ is smooth and compactly supported, $F$ is $C^k$, and and $\F$ is $C^{k+2}$, we can easily obtain the bounds
$$ \| \nu^2 \Delta_y \phi^{(0)} \|_{H^k} \leq C \nu^2 (1 + |s|)^k$$
and
$$ \| \phi^{(0)} \|_{H^k} + \| \phi^{(0)} \|_{C^k} \leq C (1 + |s|)^k.$$
Thus, arguing as in Lemma \ref{energylemma}, we have
$$ \| F(\phi^{(0)}+w)(s) - F(\phi^{(0)})(s)\|_{H^k} 
\leq C \| w(s) \|_{H^k} ((1 + |s|)^C + \| w(s) \|_{H^k})^{p-1},$$
which by \eqref{w-bound} gives
$$ \| -\nu^2 \Delta_y \phi^{(0)} + F(\phi^{(0)}+w) - F(\phi^{(0)}) \|_{H^k}
\leq C (1 + |s|)^C ( \nu^2 + e(s) + e(s)^p ),$$
which by \eqref{energy-ineq} gives the differential inequality
$$ \partial_s e(s) \leq C (1 + |s|)^C (\nu^2 + e(s) + e(s)^p).$$
Since $e(0) = 0$, we can argue exactly as in Lemma \ref{energylemma} to obtain
$$ e(s) \leq C \nu^{3/2}$$
(for instance) for $|s| \leq c |\log \nu|^c$, and the claim then follows from \eqref{w-bound} if $\nu$ is sufficiently small.
\end{proof}

\section{NLW energy transfer when $0 < s < s_c$}\label{NLW-bignorm-sec}

We now prove Theorem \ref{thm:bignorm-wave} (and hence Theorem~\ref{thm:NLW}
in the regime $ 0< s < s_c$.  
As in the Schr\"odinger case, we fix some compactly supported smooth function $w$, 
which we can choose to be the square of another such function.  For $a \in [1/2,2]$ and $0 < \nu \ll 1$, we use Lemma \ref{energylemma-wave} to construct solutions $\phi(s,y) = \phi^{(a,\nu)}(s,y)$ to \eqref{smalldispersion-wave} with initial data $\phi(0,y) := a w(y)$.  In particular we have
\begin{equation} \label{basicnlw}
\| \phi^{(a,\nu)}(s) - \phi^{(a,0)}(s) \|_{H^k} +
\| \phi^{(a,\nu)}_s(s) - \phi^{(a,0)}_s(s) \|_{H^k} \leq C \nu
\end{equation}
where
\be{phi-def-wave}
\phi^{(a,0)}(s,y) := a \F(a^{(p-1)/2} |w(y)|^{(p-1)/2} s) w(y).
\end{equation}
Applying \eqref{wave-transform}  gives then
solutions $u(t,x) = u^{(a,\nu,\lambda)}(t,x)$ to \eqref{NLW} defined by
\be{u-wave-def}
u^{(a,\nu,\lambda)}(t,x) :=
\lambda^{-2/(p-1)} \phi^{(a,\nu)}(\lambda\rp{t}, \lambda\rp{\nu x}).
\end{equation}
In particular, we have the initial data
\be{u-initial}
u^{(a,\nu,\lambda)}(0,x) = 
a \lambda^{-2/(p-1)} w(\lambda\rp{\nu x}); \quad u^{(a,\nu,\lambda)}_t(0,x) = 0.
\end{equation}

Assume that $0 < \lambda \leq \nu \ll 1$.  Arguing exactly as in Section \ref{section:positives} we have
$$ \| u^{(a,\nu,\lambda)}(0) \|_{H^s} \leq C \lambda^{s_c - s} \nu^{s-d/2}
$$
and so once again we define $\lambda$ in terms of $\nu$ by the formula \eqref{lambda-fix}.

Now we repeat the argument in Section \ref{section:positives}.  From \eqref{phi-def-wave} we have 
$$ \partial_x^j \phi^{(a,0)}(t,x) := a w(x) t^j (a^{(p-1)/2} \nabla_x |w(x)|^{(p-1)/2})^j \F^{(j)}(a^{(p-1)/2} t |w(x)|^{p-1/2}) + O(t^{j-1})$$
for $j = 0, 1, \ldots, k$.  Since $\F$ and its derivatives only vanish on a countable set (as can be easily seen from e.g.\ the Picard uniqueness theorem for ODE) we thus have
$$ \| \phi^{(a,0)}(t) \|_{H^j} \sim t^j.$$
The rest of the argument then proceeds exactly as in Section \ref{section:positives}, 
establishing Theorem \ref{thm:bignorm-wave}  for Sobolev 
parameters $0 < s < s_c$.  
\qed

\medskip
The scaling and decoherence argument in the final two paragraphs of 
\S\ref{section:positives}
can be repeated to establish Theorem \ref{thm:NLW} for 
all  $s< s_c$.
See also Lemma~\ref{decoherence} below.

\section{NLW energy transfer when $s \le -d/2$}\label{NLW-verylows-sec}

The proof of Theorem~\ref{thm:bignorm-wave} for $s \le -d/2$  is 
almost identical to the proof
of the corresponding part of Theorem~\ref{thm:bignorm} in Section~\ref{section:verylows}.  
We claim that for any prescribed $\kappa<\infty$,
the function $w$ in \eqref{phi-def-wave}
may be chosen so that  $\widehat{w}(\xi)=O(|\xi|^\kappa)$ as $|\xi|\to 0$, 
even though this is inconsistent with our earlier insistence that $w$ be the square of
a real-valued function, without disturbing the validity of \eqref{basicnlw}.  Indeed,
solutions $\phi(s,y)$ of $-\p^2_s\phi + \nu^2\Delta_y \phi =\omega|\phi|^{p-1}\phi$
exhibit uniformly finite speed of propagation so long as $|\nu|\le 1$.
Moreover, whenever $\phi$ is a solution, then so is $-\phi$.
By taking $w$ in \eqref{phi-def-wave} to be an appropriate linear combination
of squares of $C^\infty_0$ functions whose supports are widely spaced, 
we may ensure that $\widehat{w}$ vanishes to any prescribed order at $0$.
This may be done so that the supports of the summands are sufficiently widely
spaced that the uniformly bounded speed of propagation guarantees that
the solution of $-\p^2_s\phi + \nu^2\Delta_y \phi =\omega|\phi|^{p-1}\phi$
is simply the sum of the corresponding solutions for the summands.

We claim further that, as in \S\ref{section:negatives},
even when $\widehat{w}$ vanishes to high order at the origin, 
the function $\phi^{(a,0)}(s,y)$ defined in \eqref{phi-def-wave} 
need not satisfy $[\phi^{(a,0)}(s)]\widehat{\phantom{A}}(0)=0$.
Indeed, 
$$\partial_{ss} \phi^{(a,0)}(0,y) = - a^p |w(y)|^{p-1} w(y)$$ 
will not have zero integral, with respect to $y$, for generic choices of $w$.  
This implies the claim by a simple argument which is left to the reader.

At time $0$, $u^{(a,\nu,\lambda)}$ is supported where $|x|\le C\lambda \nu\rp
\le C$. Therefore
the finite speed of propagation implies that the $H^k$ norm of $u^{(a,\nu,\lambda)}$
majorizes its $H^{k,k}$ counterpart.  The rest of the argument now proceeds 
exactly as in Section~\ref{section:verylows}.
The details are left to the reader.
\qed

\section{NLW ill-posedness in the energy space}\label{energy-sec}

We now prove Theorem \ref{thm:energy-bad}.  Fix $d \geq 3$; since $s_c > 1$, we have $p > 1 + \frac{4}{d-2}$.
Now the equation is supercritical in $H^1$, so
the worst component of the $X$ norm is in fact the term $\|u\|_{p+1}^{p+1}$, rather
than $\|\nabla u\|_{L^2}^2$.  By Sobolev embedding, the $L^{p+1}$ norm scales in the same way
as the homogeneous Sobolev norm $\dot H^s$, where
\be{s-def}
s := \frac{d}{2} - \frac{d}{p+1}.
\end{equation}
Observe that $1 < s < s_c$; indeed, we have the identity $s_c - s = \frac{2}{p+1}(s_c - 1)$, which incidentally is also the identity which ensures that the kinetic and potential components of the conserved energy \eqref{energy-scale} have the scaling.  

We shall choose $w$ to be the square of a $C^\infty_0$ function which
is supported  on the ball $|x| \leq 1$ and which equals 1 
where $|x| \leq 1/2$.
We now run the construction of Section \ref{NLW-bignorm-sec}, with $s$ defined as in \eqref{s-def}.  Thus for any $0 < \nu \ll 1$ and $0 < \eps \ll 1$, $\lambda$ is defined in terms of
$\nu$ by \eqref{lambda-fix}, which in the present situation reads
\be{lambda-fix-energy}
\lambda^{\frac{2}{p+1}(s_c - 1)} \nu^{s-d/2} = \eps.
\end{equation}
The formula \eqref{u-wave-def}
defines a family $u^{(a,\nu,\lambda)}$ of solutions of \eqref{NLW}. 

As shown in Section~\ref{NLW-bignorm-sec}, there is a bound
$$ \| u^{(a,\nu,\lambda,0)}(0) \|_{H^s} \leq C \lambda^{s_c - s} \nu^{s-d/2} = C\eps.
$$
In particular, since $s>1$, we have
$$ \| \nabla u^{(a,\nu,\lambda,0)}(0) \|_{L^2} \leq C \eps.$$
Also, \eqref{u-initial} gives
$$ \| u^{(a,\nu,\lambda,0)}_t(0) \|_{L^2} = 0$$
and
$$ \| u^{(a,\nu,\lambda,0)}(0) \|_{p+1} =
a \lambda^{-2/(p-1)} (\lambda/\nu)^{d/(p+1)} \| w \|_{p+1},$$
so by \eqref{s-def} and \eqref{lambda-fix}
$$ \| u^{(a,\nu,\lambda,0)}(0) \|_{p+1} \leq C \lambda^{-2/(p-1)} (\lambda/\nu)^{\frac{d}{2} - s} = C \eps.$$
Thus we have
$$ \| u^{(a,\nu,\lambda,0)}(0) \|_{p+1} \leq C \eps.$$
A similar argument gives
$$ \| u^{(a,\nu,\lambda,0)}(0) - u^{(a',\nu,\lambda,0)}(0) \|_{p+1} \leq C \eps |a-a'|.$$
Now we investigate the behaviour for later times.  The key lemma is the following decoherence property.

\begin{lemma}\label{decoherence}  
Suppose $|a-a'| \ll 1$.  Then there exists a time $t \sim 1/|a-a'|$ such that
$$ \| \phi^{(a,0)}(t) - \phi^{(a',0)}(t) \|_{L^{p+1}} \geq c.$$
\end{lemma}

\begin{proof}
Since $w(x) = 1$ on the ball $|x| \leq 1/2$, we have from \eqref{phi-def-wave} that
$$
\| \phi^{(a,0)}(t) - \phi^{(a',0)}(t) \|_{L^{p+1}} \geq c
|a \F(a^{(p-1)/2} t) - a' \F((a')^{(p-1)/2} t)|.$$
Since $|a - a'| \ll 1$, it thus suffices to find $t \sim 1/|a-a'|$ such that
$$ |\F(a^{(p-1)/2} t) - \F((a')^{(p-1)/2} t)| \geq c.$$
The function $\F$ is non-constant and periodic, with some period $\sim 1$.  
Thus the functions $\F(a^{(p-1)/2} t)$ and $\F((a')^{(p-1)/2} t)$ are also 
non-constant and periodic, with periods differing by $\sim |a-a'|$.  
The claim then easily follows.
\end{proof}

Combining this lemma with \eqref{u-wave-def} we obtain
$$ \| u^{(a,\nu,\lambda,0)}(t\lambda)
- u^{(a',\nu,\lambda,0)}(t\lambda) \|_{p+1} \geq c \lambda^{-2/(p-1)} (\lambda/\nu)^{\frac{d}{2} - s} = c \eps.$$
Thus if we choose $|a-a'|$ sufficiently small and $\nu$ (and hence $\lambda$) sufficiently small, the claim follows.  This proves Theorem \ref{thm:energy-bad}.
\qed

\section{One-dimensional ill-posedness for NLW}\label{1d:sec}

\subsection{Well-posedness theory}
In this section we prove Theorem \ref{thm:bignorm-1d}.  
Fix $d=1$ and $\omega = \pm 1$.
The construction of solutions which demonstrate ill-posedness will rely on the
local well-posedness theory of \eqref{NLW} at the Sobolev exponent
$s_{\rm sob} = \max(0,\frac{1}{2} - \frac{1}{p})$,
so we begin with a review of that theory.
For simplicity we consider only the case of zero initial velocity $u_1 \equiv 0$.

Suppose that we have an initial data $u_0$ with norm
$$ \| u_0 \|_{H^{s_{\rm sob}}} \leq \delta$$
for some small $0 < \delta \ll 1$ to be chosen later, and we wish to solve the Cauchy problem \eqref{NLW} (with initial velocity $u_1 \equiv 0$) on the time interval $0 \leq t \leq 1$.  We shall do this by the usual iteration method, letting 
$$u^{(0)}(t,x) := \tfrac12 u_0(x+t) + \tfrac12 u_0(x-t)$$
be the solution to the free wave equation with initial data $u_0$ and initial velocity 0, and then defining iteratively for $n=1,2,\ldots$
$$ u^{(n)}(t,x) = u^{(0)}(t,x) + \tfrac{1}{2} \int_0^t \int_{x-s}^{x+s} F(u^{(n-1)}(t-s, y))\, dy\, ds,$$
where $F$ is as in \eqref{F-def}. 
Thus $u^{(n)}$ is the solution to the Cauchy problem
$$
\left\{
\begin{aligned}
\square u^{(n)}(t,x) &= F(u^{(n-1)})(t,x)
\\ 
u^{(n)}(0,x) &= u_0(x)
\\
\p_t u^{(n)}(0,x) &= 0.
\end{aligned}
\right.
$$

Observe from Sobolev embedding that $H^{s_{\rm sob}}$ controls $L^p$ if $p\ge 2$
while $L^1$ controls $H^{s_{\rm sob}-1}$ (with some regularity to spare)
since $s_{\rm sob}-1<-\tfrac12$.  
When $p<2$, $H^{s_{\rm sob}}=L^2$ controls the $L^p$ norm, provided that $x$ is
restricted to a set of uniformly bounded diameter. But we may always reduce to
this case, provided that $t\in[0,1]$, by exploiting the bounded speed of propagation
of solutions; we will freely invoke this reduction below without further mention.
Thus on the spacetime slab $[0,1] \times \R$, 
\begin{align*}
\| u^{(n)} \|_{L^\infty_t H^{s_{\rm sob}}_x}
&\leq C\| u_0 \|_{H^{s_{\rm sob}}_x}
+ C\| F(u^{(n-1)}) \|_{L^\infty_t H^{s_{\rm sob}-1}_x} \\
&\leq C\delta + C\| F(u^{(n-1)}) \|_{L^\infty_t L^1_x} \\
&\leq C\delta + C\| u^{(n-1)} \|_{L^\infty_t L^p_x}^p \\
&\leq C\delta + C\| u^{(n-1)} \|_{L^\infty_t H^{s_{\rm sob}}_x}^p.
\end{align*}
Also, we clearly have
$$ \| u^{(0)} \|_{L^\infty_t H^{s_{\rm sob}}_x} \leq C\delta.$$
Thus, if $\delta$ is sufficiently small, we obtain inductively the bound
$$ \| u^{(n)} \|_{L^\infty_t H^{s_{\rm sob}}_x} \leq C\delta$$
uniformly in $n$.  

Next, we study the decay of successive differences $u^{(n+1)} - u^{(n)}$.  We
adopt the convention that $u^{(-1)} \equiv 0$.  Adapting the above 
estimate to $u^{(n+1)} - u^{(n)}$ in the usual manner, we obtain
\begin{align*}
\| u^{(n+1)} - u^{(n)}\|_{L^\infty_t H^{s_{\rm sob}}_x}
&\leq C\| F(u^{(n)}) - F(u^{(n-1)}) \|_{L^\infty_t H^{s_{\rm sob}-1}_x} \\
&\leq C\| F(u^{(n)}) - F(u^{(n-1)})  \|_{L^\infty_t L^1_x} \\
&\leq C( \| u^{(n)} \|_{L^\infty_t L^p_x} + \| u^{(n-1)} \|_{L^\infty_t L^p_x})^{p-1} 
\| u^{(n)} - u^{(n-1)} \|_{L^\infty_t L^p_x} \\
&\leq C\delta^{p-1} \| u^{(n)} - u^{(n-1)} \|_{L^\infty_t H^{s_{\rm sob}}_x}.
\end{align*}
Thus we obtain a bound of the form
$$ \| u^{(n)} - u^{(n-1)} \|_{L^\infty_t H^{s_{\rm sob}}_x} \leq C(C\delta)^{n(p-1)} \delta.$$
In particular, the $u^{(n)}$ converge in $L^\infty_t H^{s_{\rm sob}}_x$ to a function $u(t,x)$ with the convergence bounds
\be{convergence} \| u - u^{(n)} \|_{L^\infty_t H^{s_{\rm sob}}_x} \leq C(C\delta)^{(n+1)(p-1)} \delta.
\end{equation}
It is easy to verify, via a contraction mapping argument, 
that $u$ is the unique $C^0_tH^{s_{\rm sob}}_x$ solution to \eqref{NLW}.  

These bounds can be improved.
Note that $s_{\rm sob}$ is strictly $<\tfrac12$ for any $p$.
For any $r\in(0,\tfrac12)$, and in particular for certain $r>s_{\rm sob}$,
\begin{align}
\| u^{(n+1)} - u^{(n)}\|_{L^\infty_t H^{r}_x}
&\leq C
\| F(u^{(n)}) - F(u^{(n-1)}) \|_{L^\infty_t H^{r-1}_x} 
\notag
\\
\notag
&\leq C\| F(u^{(n)}) - F(u^{(n-1)})  \|_{L^\infty_t L^1_x} \\
\notag
&\leq C\delta^{p-1} \| u^{(n)} - u^{(n-1)} \|_{L^\infty_t H^{s_{\rm sob}}_x}
\\
& \le C\delta^{(p-1)(n+1)+1}. \label{1dWPimproved} 
\end{align}

\subsection{Proof of Theorem~\ref{thm:bignorm-1d} for $p\ge 2$}
We next apply this well-posedness theory to a specific class of initial data.  
Let $w\in C^\infty_0$ be smooth and compactly supported, not identically zero,
with $\widehat{w}(\xi)= O(|\xi|^M)$ as $\xi\to 0$, where $M$ is sufficiently large
for later purposes. There exist such $w$ with the additional property that for all 
sufficiently large positive $N$,
\begin{equation} \label{bigiterate1}
N\big|
\int_0^1 \int_{x-s}^{x+s} F(w(N(y-1+s))+w(N(y+1-s)))\,dy\,ds
\big|
\ge c \text{ uniformly for all } |x|\le\tfrac12
\end{equation}
for some constant $c>0$ independent of $N$.
Indeed, if $w$ is even then elementary calculations show that the above
expression equals
\begin{equation}
2\omega\int_\reals |w(y)|^{p-1}w(y)\,dy + O(N^{-1})
\end{equation}
uniformly for $|x|\le\tfrac12$. We leave it to the reader to verify, by a simple variational
argument, that the leading term can be nonzero.
Henceforth we assume that $w$ satisfies \eqref{bigiterate1}.

Consider the initial data
\be{u0-def}
 u_0(x) := \delta N^{1/p} w(Nx); \quad u_1(x) := 0
\end{equation}
where $0 < \delta \ll 1$ is a small absolute constant (depending only on $p$, $w$) and $N \gg 1$ is a parameter to be chosen later.  A calculation gives
$$ \| u_0 \|_2 \leq C\delta N^{1/p - 1/2}$$
and
$$ \| u_0 \|_{H^1} \leq C\delta N^{1/p + 1/2}$$
so by convexity and the definition of $s_{\rm sob}$, we see that
$$ \| u_0 \|_{H^{s_{\rm sob}}} \leq C\delta.$$
Thus by the previous analysis, if $\delta$ is sufficiently small then we have a solution $u(t,x)$ to \eqref{NLW} on the spacetime slab $[0,1] \times \R$.  
Furthermore, by \eqref{convergence} 
$$\| u - u^{(1)} \|_{L^\infty_t H^{s_{\rm sob}}_x} \leq C\delta^{2p-1}, $$
and in particular, for any $s\le s_{\rm sob}$, 
\be{convergence-1}
\| u(1) - u^{(1)}(1) \|_{H^{s}_x} \leq C\delta^{2p-1}.
\end{equation}
Recall that $u^{(1)}$ is given by the formula
\bas
u^{(1)}(t,x) &= u^{(0)}(t,x) + \tfrac{1}{2} \int_0^t \int_{x-s}^{x+s} F(u^{(0)})(t-s, y)\ dy ds \\
&= \frac{u_0(x-t) + u_0(t+x)}{2} + \tfrac{1}{2} \int_0^t \int_{x-s}^{x+s} F(\frac{u_0(y-t+s) + u_0(y+t-s)}{2})\ dy ds.
\end{align*}
We now compute $u^{(1)}$ on the region $t=1$, $|x| < \tfrac12$.  
In this region $u_0(x-t)$ and $u_0(t+x)$ vanish (if $N$ is large enough), and a 
brief calculation (using \eqref{u0-def} and \eqref{bigiterate1}) reveals a lower bound
$$ |u^{(1)}(1,x)| \geq c \delta^p\ \text{ for all } |x|<\tfrac12$$
for all $N$ greater than some fixed constant independent of $\delta$.
Let 
\[U = u^{(1)}-u^{(0)}.\]
Since the supports of $u^{(1)}(1), u^{(0)}(1)$ 
are contained in a fixed compact set, independent of $N$, we have 
\[\|U(1)\|_{L^1}\le C\|U(1)\|_{L^2}=O(\delta^p),\] 
and hence by \eqref{1dWPimproved}
for any constants $0<a\le 1\le A$ and any $r<\tfrac12$,
\begin{equation}
\label{u1-large}
\int_{a<|\xi|<A} |\widehat{ U(1) }(\xi)|^2\,d\xi
\ge c\delta^{2p} - a\delta^{2p} - A^{-2r}\|U(1)\|_{H^r}^2
\ge c\delta^{2p} - Ca\delta^{2p} - CA^{-2r}\delta^{2p}
\end{equation}
uniformly in $N,\delta$ provided that $N\ge 1$ and $\delta\le 1$.
Choose $a,A$ sufficiently small and large, respectively,
so that this is $\ge c'\delta^{2p}$.
Now $\widehat{u^{(0)}(t)} = O(\delta N^{\tfrac1p-1})$ in $L^\infty$
uniformly for all $t$, and $p>1$, so the same lower bound holds for
$\widehat{u^{(1)}(1)}$ uniformly for all sufficiently large $N$.
From \eqref{convergence-1} it follows that $u(1)$ itself satisfies 
\begin{equation}
\int_{a<|\xi|<A} |\widehat{ u{(1)}}(\xi)|^2\,d\xi
\ge c'\delta^{2p}
\end{equation}
for some $c'>0$,
provided that $\delta$ is first chosen to be 
sufficiently small but fixed, and $N$ is then taken to be sufficiently large.


Now we utilize the scaling symmetry \eqref{wave-scale}, defining
$$ u^{\lambda}(t,x) := \lambda^{-2/(p-1)} u(\frac{t}{\lambda}, \frac{x}{\lambda})$$
for some parameter $0 < \lambda \le 1$ to be chosen later.  This function $u^{\lambda}$ 
satisfies the wave equation \eqref{NLW} on the slab $[0,\lambda] \times \R$, but with initial datum
$$ u^\lambda(0,x) = \lambda^{-2/(p-1)} \delta N^{1/p} w(Nx/\lambda); \quad u^\lambda_t(0,x) = 0.$$
A simple calculation then shows that
\begin{equation}
\| u^\lambda(0) \|_{H^s} \leq C\lambda^{-2/(p-1)} \delta N^{1/p} (N/\lambda)^{s-1/2}
\end{equation}
provided that $M$ is chosen so that $-M<s$, which we will assume henceforth.
Thus if $p\ge 2$,
\be{initial}
\| u^\lambda(0) \|_{H^s} \leq C\delta \lambda^{s_c - s} N^{s - s_{\rm sob}},
\end{equation}
while for $p<2$ there is the weaker bound $\le C\delta\lambda^{s_c - s} N^{s-(\frac12-\frac1p)}$.

On the other hand, \eqref{u1-large} and rescaling give us
$$ \| u^\lambda(\lambda) \|_{H^s}^2 
\ge c\lambda^{-2s}\int_{a<\lambda|\xi|<A} 
|\widehat{u^\lambda(\lambda)}(\xi)|^2
\,d\xi
\gtrsim \lambda^{-2s}\delta^{2p}
\lambda\cdot\lambda^{-4/(p-1)} $$
by \eqref{u1-large} provided that $a,A$ are chosen as above.
Hence
$$\|u^\lambda(\lambda)\|_{H^s} \ge c\delta^p \lambda^{1/2}\lambda^{-2/(p-1)}
\lambda^{-s}
= c\delta^p \lambda^{s_c-s}.
$$

We can now conclude the proof of Theorem~\ref{thm:bignorm-1d}
in the case where $s_{\rm sob}\ge 0$, or equivalently $p\ge 2$.
Suppose first that $s>s_c$.
Fixing a sufficiently small $0<\delta\ll 1$, we may choose $a,A$ as above.
We may ensure that $u^\lambda(\lambda)$ is arbitrarily large in $H^s$ norm 
by choosing $\lambda$ sufficiently small. 
If $N$ is next chosen to be sufficiently large depending on $\delta$
and $\lambda$, then \eqref{initial}  guarantees that
the initial datum becomes arbitrarily small in the $H^s\times H^{s-1}$ norm. 
All the conclusions of Theorem \ref{thm:bignorm-1d} for $s>s_c$ follow directly.

For $s\le s_c$ we choose $\lambda$ to be an arbitrarily small constant.
As $N\to\infty$, these estimates then yield the weaker conclusion announced
in Theorem~\ref{thm:bignorm-1d}:
initial data which are arbitrarily small in $H^s\times H^{s-1}$
can give rise to solutions which have $H^s$ norms comparable to a positive power of $\lambda$
in time $\lambda$.

\subsection{The case $p<2$.}
Next we prove Theorem~\ref{thm:bignorm-1d}
when $p<2$, so that
$s_{\rm sob} = 0$. The reasoning above continues to apply whenever $s<\tfrac12-\tfrac1p$, since 
\eqref{initial} holds with $N^{s-s_{\rm sob}}$ replaced by $N^{s-(\tfrac12-\tfrac1p)}$,
but another construction seems necessary when
$\tfrac12-\tfrac1p\le s<0$. We may thus assume that $-\tfrac12<s$.
Again the local well-posedness theory 
for $s=s_{\rm sob}$ will be used, but the initial data $u_0$ (and hence the solution $u$)
will be modified.

We shall need a smooth function $\phi: \R \to \R$ which is periodic of period 1, 
which has the moment property
$$ \int_0^1 \phi(v)\ dv = 0 $$
but also satisfies
\begin{equation}\label{f-nonzero}
\int_0^1 \int_0^1 F(\phi(v) + \phi(w))\ dv dw = c_0 \neq 0
\end{equation}
for non-zero some constant $c_0$, where $F(z) = |z|^{p-1} z$.  Such functions exist;
for instance we can set 
$$\phi(v) := A \Phi(Av) - 1 \hbox{ for } v \in [0,1]$$
and extended periodically for all $v$, where $A \gg 1$ is a sufficiently large parameter, and $\Phi$ is a fixed non-negative bump function on $[0,1]$ with total mass $1$.

We begin again with initial data of the form
$$ u_0(x) := \delta \phi(Nx) \psi(x); \quad u_1(x) := 0$$
where $N \gg 1$ is a large parameter, $0 < \delta \ll 1$ is a small parameter, and $\psi$ is a 
non-negative smooth function supported on the interval $[-4,4]$ which equals 1 on $[-2,2]$.

Let us first compute the $H^s$ norm of $u_0$.  Clearly $\|u_0\|_2 = O(\delta)$; since $s_{\rm sob} = 0$, this means that the local well-posedness theory of the previous section will apply in $L^2$ 
if $\delta$ is small enough.  Also, since $\phi$ has mean zero, we may write $\phi$ as a derivative $\phi = \Psi_x$ for some smooth periodic function $\Psi$.  Thus
$$ u_0 = \frac{\delta}{N} \frac{d}{dx} (\Psi(Nx) \psi(x)) - \frac{\delta}{N} \Psi(Nx) \psi_x(x).$$
The second term has an $L^2$ norm of $O(\delta/N)$, and hence an $H^s$ norm of $O(\delta/N)$.  The first term has an $L^2$ norm of $O(\delta)$, and an $\dot H^{-1}$ norm of $O(\delta/N)$, and hence has an $H^s$ norm of $O(\delta N^s)$ since $-1 <s < 0$.  Putting this all together we see that
$$ \| u_0 \|_{H^s} \leq C \delta N^s,$$
and so the $H^s$ norm can be made small by choosing $N$ sufficiently large depending on $\delta$.

It is necessary to analyze
some iterates of $u$.  From energy estimates (or just by translation invariance) we have
$$ \| u^{(0)}(t) \|_{H^s} \leq C \| u_0 \|_{H^s} \leq C \delta N^s$$
for all times $t$.  Next consider $u^{(1)}$.  As before, 
$$ u^{(1)}(t,x) =
u^{(0)}(t,x) +
\tfrac{1}{2} \int_0^t \int_{x-s}^{x+s}
F\big(\tfrac12{u_0(y-t+s) + \tfrac12u_0(y+t-s)}\big) dy ds.$$
Let us restrict attention to the region $t=1$, $|x| < 1/2$.  Then $\psi$ can be replaced by 1 throughout.  Expanding the second term, and exploiting the homogeneity of $F$, we thus have
$$ u^{(1)}(1,x)
= u^{(0)}(1,x) + C \delta^p
\int_0^1 \int_{x-s}^{x+s}
F(\phi(N(y-1+s)) + \phi(N(y+1-s)))\ dy ds$$
where $C$ is a positive absolute constant (depending only on $p$ and $\phi$).  Making the change of variables $v := N(y-1+s)$, $w := N(y+1-s)$, this becomes
$$ u^{(1)}(1,x) = u^{(0)}(1,x) +
C \delta^p N^{-2} \iint_{N(x-1) \leq v \leq w \leq N(x+1)} F(\phi(v) + \phi(w))\,dv\,dw,$$
where $C$ again denotes a positive constant.  Using \eqref{f-nonzero} and the periodicity, 
partitioning the region of integration into unit squares plus a remaining region of area $O(N)$, 
we thus have
\begin{equation} \label{p<2control}
 u^{(1)}(1,x) - u^{(0)}(1,x) =
c \delta^p + O(\delta^p N^{-1})
\end{equation}
for some positive $c$.

Since $U(1)$ is supported in a fixed bounded set
independent of $\delta,N$, we have 
\begin{equation}
\int_{a<|\xi|<A} |\widehat{U(1)}(\xi)|^2\,d\xi
\ge c\delta^{2p}
\end{equation}
by \eqref{p<2control} and
the same reasoning as given above for the case $p\ge 2$,
provided that $a,A$ are chosen to be sufficiently small and large,
respectively.
On the other hand, 
$|\widehat{u^{(0)}(1)}(\xi)|
\le C_A\delta N^{-1}$  for all $|\xi|\le A$, uniformly
for all $0<\delta\le 1$ and all $N\ge 1$.
By combining these bounds with \eqref{p<2control}
 we find that for any sufficiently small $\delta>0$, 
for all sufficiently large $N$,
\begin{equation}
\int_{a<|\xi|<A} |\widehat{u^{(1)}(1)}(\xi)|^2\,d\xi
\ge c\delta^{2p}.
\end{equation}

With this lower bound in hand
we can conclude the proof of Theorem~\ref{thm:bignorm-1d} in the
case $p< 2$ by repeating the final part
of the argument given above for $p\ge 2$,
using the local well-posedness theory to approximate $u$ by $u^{(1)}$, 
then scaling to produce a richer family of solutions.

\subsection{Proof of Corollary~\ref{1dcorollary}}
The corollary is proved by considering
initial data of product form $\eta(x')f(x_d)$ where $x=(x',x_d)\in\reals^{d-1}\times\reals$,
$f$ is as in the proof of Theorem~\ref{thm:bignorm-1d}, and $\eta$ is a fixed $C^\infty$,
compactly supported function which is $\equiv 1$ on a sufficiently large ball in
$\reals^{d-1}$. By finite speed of propagation, the corresponding solutions, 
assuming existence and uniqueness, will likewise have product form for $x'$ in a 
fixed smaller ball, so norm inflation in $\reals^d$ follows from the growth
already established in $\reals^1$. The details of this argument are left to the reader.
It may be possible to obtain more refined results by allowing the support of $\eta$
to shrink with that of $f$, but we have not investigated this.

\end{document}